\documentclass[12pt]{amsart}
\usepackage{amssymb,amsmath,amsthm,latexsym,booktabs,todonotes,color,comment, eucal,uncial,eufrak,pict2e,url}
\usepackage[normalem]{ulem}
\newcommand{\stkout}[1]{\ifmmode\text{\sout{\ensuremath{#1}}}\else\sout{#1}\fi}
\theoremstyle{definition}
\newtheorem{definition}{Definition}
\newtheorem{notation}[definition]{Notation}

\newtheorem{remark}[definition]{Remark}
\newtheorem{problem}[definition]{Problem}

\theoremstyle{plain}
\newtheorem{lemma}[definition]{Lemma}
\newtheorem{proposition}[definition]{Proposition}
\newtheorem{theorem}[definition]{Theorem}

\newtheorem{conjecture}[definition]{Conjecture}
\newtheorem{example}[definition]{Example}

\begin{document}
\title[Odd order $C_4$-face-magic projective grid graphs]{Odd order $C_4$-face-magic $m \times n$ projective grid graphs having $C_4$-face-magic value $2mn+1$ or $2mn+3$}
\author{Stephen J. Curran}
\address{Department of Mathematics\\
         University of Pittsburgh at Johnstown\\
         Johnstown, PA 15904\\
         USA}
         \email{sjcurran@pitt.edu}

\date{June 7, 2022.  }

\begin{abstract}
For a graph $G = (V, E)$ embedded in the projective plane, let $\mathcal{F}(G)$ denote the set of faces of $G$.
Then, $G$ is called a $C_n$-face-magic projective graph if there exists a bijection $f: V(G) \to \{1, 2, \dots, |V(G)|\}$
such that for any $F \in \mathcal{F}(G)$ with $F \cong C_n$, the sum of all the vertex labels around $C_n$ is a constant $S$.
We consider the $m \times n$ grid graph, denoted by $\mathcal{P}_{m,n}$, embedded in the projective plane in the natural way.

Let $m \geqslant 3$ and $n \geqslant 3$ be odd integers.
It is known that the $C_4$-face-magic value of a $C_4$-face-magic labeling on $\mathcal{P}_{m,n}$ is either $2mn+1$, $2mn+2$, or $2mn+3$.
The characterization of $C_4$-face-magic labelings on $\mathcal{P}_{m,n}$ having $C_4$-face-magic value $2mn+2$ is known.
In this paper, we determine a category of $C_4$-face-magic labelings on $\mathcal{P}_{m,n}$ for which the
$C_4$-face-magic value is either $2mn+1$ or $2mn+3$.
It is conjectured that these are the only $C_4$-face-magic labeling on $\mathcal{P}_{m,n}$ having $C_4$-face-magic value $2mn+1$ or $2mn+3$.
\end{abstract}

\maketitle

\section{Introduction} \label{S:introduction}

Graph labelings were introduced by Kotzig and Rosa \cite{Kotzig} in the 1970s.
Various applications to graph labelings include graph decomposition problems, radar pulse code designs, X-ray
crystallography and communication network models.
We refer the interested reader to J. A. Gallian's comprehensive dynamic survey on graph labelings \cite{Gallian}
for further investigation.

The reader should consult Chartrand, Lesniak, and Zhang \cite{Chartrand} for concepts and notation not explicitly defined in this paper.
The graphs in this paper are connected multigraphs.
The concept of a $C_4$-face-magic labeling was first applied to planar graphs.
We apply the concept to graphs embedded on a projective plane.
For a planar (toroidal, Klein bottle, projective) graph $G = (V, E)$ embedded in the plane (torus, Klein bottle, projective plane),
let $\mathcal{F}(G)$ denote the set of faces of $G$.
Then, $G$ is called a \textit{$C_n$-face-magic planar (toroidal, Klein bottle, projective)} graph if there
exists a bijection $f: V(G) \to \{1, 2, \dots, |V(G)|\}$ such that for any $F \in \mathcal{F}(G)$ with $F \cong C_n$,
the sum of all the vertex labels aroiund $C_n$ is a constant $S$. Here, the constant $S$ is called
a \textit{$C_n$-face-magic value} of $G$.
More generally, $C_n$-face-magic planar graph labelings are a special case of $(a, b, c)$-magic labeling
introduced by Lih \cite{Lih}.
For assorted values of $a, b$ and $c$, Baca and others \cite{Baca1, Baca2, Baca3, Hsieh, Kasif, Kath, Lih} have analyzed
the problem for various classes of graphs.
Wang \cite{Wang} showed that the toroidal grid graphs $C_m \times C_n$ are antimagic for all integers $m,n\geqslant 3$.
Butt et al. \cite{Buttetal} investigated face antimagic labelings on toroidal and Klein bottle grid graphs.

Curran, Low and Locke \cite{CurranLow, CurranLowLocke1} investigated $C_4$-face-magic labelings
on an $m \times n$ toroidal grid graph.
They showed that $C_m \times C_n$ admits a $C_4$-face-magic labeling if and only if either $m=2$, or $n=2$, or both $m$ and $n$ are even.
Curran, Low and Locke \cite{CurranLowLocke2} also examined $C_4$-face-magic labelings on an $m \times n$
Klein bottle grid graph.
They showed that an $m \times n$ Klein bottle grid graph  admits a $C_4$-face-magic labeling if and only if $n$ is even.
In this paper, we consider $C_4$-face-magic projective labelings on an $m \times n$
projective grid graph.

Curran \cite{Curran2} showed that an $m \times n$ projective grid graph  admits a $C_4$-face-magic labeling if and only if both $m$ and $n$ have the same parity.
When $m$ and $n$ are even, the $C_4$-face-magic value of a $C_4$-face-magic labeling on an $m\times n$ projective grid graph must be $2mn+2$.
Also, when $m$ and $n$ are odd, the $C_4$-face-magic value of a $C_4$-face-magic labeling on an $m\times n$ projective grid graph is
either $2mn+1$, $2mn+2$, or $2mn+3$.
The $C_4$-face-magic labelings on $\mathcal{P}_{m,n}$ having $C_4$-face-magic value $2mn+2$ were characterized in \cite{Curran2}.
In this paper, we determine a category of the $C_4$-face-magic labelings on $\mathcal{P}_{m,n}$ having $C_4$-face-magic value $2mn+1$ or $2mn+3$.
We conjecture that these are the only $C_4$-face-magic labelings on $\mathcal{P}_{m,n}$ having  $C_4$-face-magic value $2mn+1$ or $2mn+3$.

\section{Preliminaries} \label{S:preliminaries}

In this section we introduce definitions and known results about $C_4$-face-magic projective labelings on an $m \times n$ projective grid graph.

\begin{definition}
Let $m$ and $n$ be integers such that $m,n\geqslant 2$.
The {\it $m \times n$ projective grid graph}, denoted by $\mathcal{P}_{m,n}$, is the graph whose vertex set is
\begin{align*}
V\left( \mathcal{P}_{m,n} \right)
=\left\{ \left( i,j\right) :1\leqslant i\leqslant m,1\leqslant j\leqslant n\right\},
\end{align*}
and whose edge set consists of the following edges:
\begin{itemize}
\item there is an edge from $(i,j)$ to $(i,j+1)$, for $1\leqslant i \leqslant m$ and $1\leqslant j \leqslant n-1$,
\item there is an edge from $(i,n)$ to $(m+1-i,1)$, for $1\leqslant i \leqslant m$,
\item there is an edge from $(i,j)$ to $(i+1,j)$, for $1\leqslant i \leqslant m-1$ and $1\leqslant j \leqslant n$, and
\item there is an edge from $(m,j)$ to $(1,n+1-j)$, for $1\leqslant j \leqslant n$.
\end{itemize}
The graph $\mathcal{P}_{m,n}$ has a natural embedding on the projective plane.
This graph is a multigraph since there are double edges on the vertex sets $\{ (1,1), (m,n)\}$ and $\{ (m,1), (1,n)\}$.
\end{definition}

\begin{example} \label{exampleFiveByFiveProjectiveGridGraph}
  The $5 \times 5$ projective grid graph $\mathcal{P}_{5,5}$ is illustrated in Fig. \ref{figFiveByFiveProjectiveGridGraph}.
  Due to the orientation of the vertices in $\mathcal{P}_{m,n}$, we refer to the vertices
  $\{ (i,j): 1 \leqslant j \leqslant n\}$ as column $i$ of $V(\mathcal{P}_{m,n})$ and
  $\{ (i,j): 1 \leqslant i \leqslant m\}$ as row $j$ of $V(\mathcal{P}_{m,n})$.
  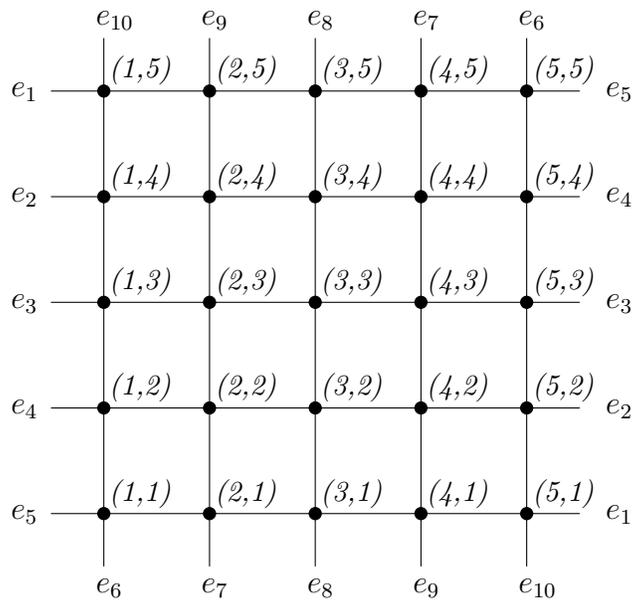
\begin{figure}
\begin{picture}(225,225)(-25,-25)

\multiput(0,0)(40,0){5}{\circle*{5}}
\multiput(0,40)(40,0){5}{\circle*{5}}
\multiput(0,80)(40,0){5}{\circle*{5}}
\multiput(0,120)(40,0){5}{\circle*{5}}
\multiput(0,160)(40,0){5}{\circle*{5}}

\put(-20,0){\line(1,0){200}}
\put(-20,40){\line(1,0){200}}
\put(-20,80){\line(1,0){200}}
\put(-20,120){\line(1,0){200}}
\put(-20,160){\line(1,0){200}}

\put(0,-20){\line(0,1){200}}
\put(40,-20){\line(0,1){200}}
\put(80,-20){\line(0,1){200}}
\put(120,-20){\line(0,1){200}}
\put(160,-20){\line(0,1){200}}

\put(2,5){\small (1,1)}
\put(42,5){\small (2,1)}
\put(82,5){\small (3,1)}
\put(122,5){\small (4,1)}
\put(162,5){\small (5,1)}

\put(2,45){\small (1,2)}
\put(42,45){\small (2,2)}
\put(82,45){\small (3,2)}
\put(122,45){\small (4,2)}
\put(162,45){\small (5,2)}

\put(2,85){\small (1,3)}
\put(42,85){\small (2,3)}
\put(82,85){\small (3,3)}
\put(122,85){\small (4,3)}
\put(162,85){\small (5,3)}

\put(2,125){\small (1,4)}
\put(42,125){\small (2,4)}
\put(82,125){\small (3,4)}
\put(122,125){\small (4,4)}
\put(162,125){\small (5,4)}

\put(2,165){\small (1,5)}
\put(42,165){\small (2,5)}
\put(82,165){\small (3,5)}
\put(122,165){\small (4,5)}
\put(162,165){\small (5,5)}

\put(-35,158){$e_{1}$}
\put(-35,118){$e_{2}$}
\put(-35,78){$e_{3}$}
\put(-35,38){$e_{4}$}
\put(-35,-2){$e_{5}$}

\put(-3,-30){$e_{6}$}
\put(37,-30){$e_{7}$}
\put(77,-30){$e_{8}$}
\put(117,-30){$e_{9}$}
\put(157,-30){$e_{10}$}

\put(190,158){$e_{5}$}
\put(190,118){$e_{4}$}
\put(190,78){$e_{3}$}
\put(190,38){$e_{2}$}
\put(190,-2){$e_{1}$}

\put(-3,185){$e_{10}$}
\put(37,185){$e_{9}$}
\put(77,185){$e_{8}$}
\put(117,185){$e_{7}$}
\put(157,185){$e_{6}$}
\end{picture}
\caption{$5 \times 5$ projective grid graph $\mathcal{P}_{5,5}$.}
\label{figFiveByFiveProjectiveGridGraph}
\end{figure}

\end{example}

The following result determines when $\mathcal{P}_{m,n}$ admits a $C_4$-face-magic projective labeling.

\begin{theorem}[\cite{Curran2}, Theorem 10] \label{thmProjectiveGridIsC4FaceMaigc}
Let $m$ and $n$ be integers such that $m,n\geqslant 2$. Then $\mathcal{P}_{m,n}$ admits
a $C_4$-face-magic projective labeling if and only if $m$ and $n$ have the same parity.
\end{theorem}

Furthermore, one can determine the possible $C_4$-face-magic values of a $C_4$-face-magic projective labeling on $\mathcal{P}_{m,n}$.
The next lemma determines the $C_4$-face-magic value of a $C_4$-face-magic projective labeling on $\mathcal{P}_{m,n}$ when $m$ and $n$ are even integers.

\begin{lemma}[\cite{Curran2}, Lemma 5] \label{lemmaEvenByEvenProjectiveGridC4FaceValue}
Suppose $m \geqslant 2$ and $n \geqslant 2$ are even integers.
Let $\{ x_{i,j} : (i,j)\in V(\mathcal{P}_{m,n}) \}$ be a $C_4$-face-magic projective labeling on $\mathcal{P}_{m,n}$ with $C_4$-face-magic value $S$.
Then $S= 2mn+2$.
\end{lemma}

In fact, one can determine the digon face sum values on the digon vertex sets $\{ (1,1), (m,n)\}$ and $\{(m1), (1,n)\}$
on $\mathcal{P}_{m,n}$ for odd integers $m$ and $n$.

\begin{lemma}[\cite{Curran2}, Lemma 6] \label{lemmaOddByOddProjectiveGridC4FaceValue}
Let $m\geqslant 3$ and $n\geqslant 3$ be odd integers.
Let $\{ x_{i,j} : (i,j)\in V(\mathcal{P}_{m,n}) \}$ be a $C_4$-face-magic projective labeling on $\mathcal{P}_{m,n}$ with $C_4$-face-magic value $S$
Let $D_1=x_{1,1} + x_{m,n}$ and $D_2=x_{m,1} + x_{1,n}$ be the face sums of the two digons constructed from the pair of vertices at opposite corners of $\mathcal{P}_{m,n}$.
Then either
\begin{enumerate}
  \item $S=2mn+1$ and $D_1=D_2= \frac{3}{2} mn+ \frac{1}{2}$,
  \item $S=2mn+2$ and $D_1=D_2=mn+1$, or
  \item $S=2mn+3$ and $D_1=D_2= \frac{1}{2} mn+ \frac{3}{2}$.
\end{enumerate}
\end{lemma}

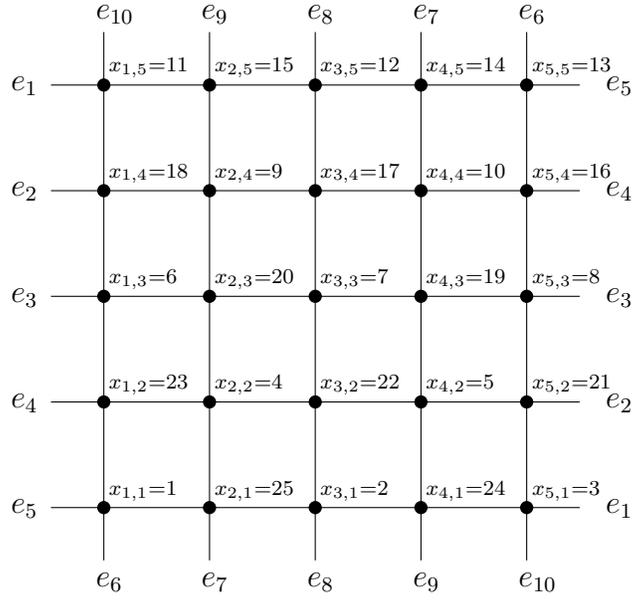
\begin{figure}
\begin{picture}(225,225)(-25,-25)

\multiput(0,0)(40,0){5}{\circle*{5}}
\multiput(0,40)(40,0){5}{\circle*{5}}
\multiput(0,80)(40,0){5}{\circle*{5}}
\multiput(0,120)(40,0){5}{\circle*{5}}
\multiput(0,160)(40,0){5}{\circle*{5}}

\put(-20,0){\line(1,0){200}}
\put(-20,40){\line(1,0){200}}
\put(-20,80){\line(1,0){200}}
\put(-20,120){\line(1,0){200}}
\put(-20,160){\line(1,0){200}}

\put(0,-20){\line(0,1){200}}
\put(40,-20){\line(0,1){200}}
\put(80,-20){\line(0,1){200}}
\put(120,-20){\line(0,1){200}}
\put(160,-20){\line(0,1){200}}

\put(2,5){${\scriptstyle x_{1,1}=1}$}
\put(42,5){${\scriptstyle x_{2,1}=25}$}
\put(82,5){${\scriptstyle x_{3,1}=2}$}
\put(122,5){${\scriptstyle x_{4,1}=24}$}
\put(162,5){${\scriptstyle x_{5,1}=3}$}

\put(2,45){${\scriptstyle x_{1,2}=23}$}
\put(42,45){${\scriptstyle x_{2,2}=4}$}
\put(82,45){${\scriptstyle x_{3,2}=22}$}
\put(122,45){${\scriptstyle x_{4,2}=5}$}
\put(162,45){${\scriptstyle x_{5,2}=21}$}

\put(2,85){${\scriptstyle x_{1,3}=6}$}
\put(42,85){${\scriptstyle x_{2,3}=20}$}
\put(82,85){${\scriptstyle x_{3,3}=7}$}
\put(122,85){${\scriptstyle x_{4,3}=19}$}
\put(162,85){${\scriptstyle x_{5,3}=8}$}

\put(2,125){${\scriptstyle x_{1,4}=18}$}
\put(42,125){${\scriptstyle x_{2,4}=9}$}
\put(82,125){${\scriptstyle x_{3,4}=17}$}
\put(122,125){${\scriptstyle x_{4,4}=10}$}
\put(162,125){${\scriptstyle x_{5,4}=16}$}

\put(2,165){${\scriptstyle x_{1,5}=11}$}
\put(42,165){${\scriptstyle x_{2,5}=15}$}
\put(82,165){${\scriptstyle x_{3,5}=12}$}
\put(122,165){${\scriptstyle x_{4,5}=14}$}
\put(162,165){${\scriptstyle x_{5,5}=13}$}

\put(-35,158){$e_{1}$}
\put(-35,118){$e_{2}$}
\put(-35,78){$e_{3}$}
\put(-35,38){$e_{4}$}
\put(-35,-2){$e_{5}$}

\put(-3,-30){$e_{6}$}
\put(37,-30){$e_{7}$}
\put(77,-30){$e_{8}$}
\put(117,-30){$e_{9}$}
\put(157,-30){$e_{10}$}

\put(190,158){$e_{5}$}
\put(190,118){$e_{4}$}
\put(190,78){$e_{3}$}
\put(190,38){$e_{2}$}
\put(190,-2){$e_{1}$}

\put(-3,185){$e_{10}$}
\put(37,185){$e_{9}$}
\put(77,185){$e_{8}$}
\put(117,185){$e_{7}$}
\put(157,185){$e_{6}$}
\end{picture}
\caption{$C_4$-face-magic projective labeling on $\mathcal{P}_{5,5}$ having $C_4$-face-magic value $53$.}
\label{figFiveByFiveProjectiveGridGraphLabeling}
\end{figure}

\begin{example}
Fig. \ref{figFiveByFiveProjectiveGridGraphLabeling} illustrates a $C_4$-face-magic projective labeling on the $5 \times 5$ projective grid graph $\mathcal{P}_{5,5}$
that has $C_4$-face-magic value 53.
\end{example}

For odd integers $m \geqslant 3$ and $n \geqslant 3$, the $C_4$-face-magic projective labelings on $\mathcal{P}_{m,n}$ having $C_4$-face-magic
value $2mn+2$ were characterized in \cite{Curran2}. The statement of this characterization involves several technical definitions.
So, we refer the reader to \cite{Curran2} for the details of this characterization.
However, it is relatively easy to count the number of distinct $C_4$-face-magic projective
labelings on $\mathcal{P}_{m,n}$ having $C_4$-face-magic value $2mn+2$.
We need the following definition in order to sate this result.

\begin{definition} \label{defnProjectiveFactorizationSequence}
  Suppose there exists a positive integer $k$ such that one of the two following conditions hold.
  \begin{enumerate}
    \item There are factorizations of $m=m_1 m_2 \ldots m_{k}$ and $n=n_1 n_2 \ldots n_k$, where $m_i >1$ and $n_i>1$ for all $1 \leqslant i \leqslant k$.
    \item There are factorizations of $m=m'_1 m'_2 \ldots m'_{k}m'_{k+1}$ and $n=n'_1 n'_2 \ldots n'_{k}$,
    where $m'_i >1$ for all $1 \leqslant i \leqslant k+1$, and $n'_i>1$ for all $1 \leqslant i \leqslant k$.
  \end{enumerate}
  We say that $(m_1,n_1,m_2,n_2,\ldots,m_k,n_k)$ is an $(m,n)$-\textit{projective factorization sequence of length $2k$}.
  Also, we say $(m'_1,n'_1,m'_2,n'_2,\ldots,m'_{k},n'_{k},m'_{k+1})$ is an $(m,n)$-\textit{projective factorization sequence of length $2k+1$}.
  For convenience, we let $n'_{k+1}=1$ and refer to
  $(m'_1,n'_1,m'_2,\allowbreak n'_2,\allowbreak \ldots,\allowbreak m'_{k+1},\allowbreak n'_{k+1})$ as an $(m,n)$-\textit{projective factorization sequence of length $2k+1$}.
  In addition, we say that $(m_1,n_1,m_2,n_2,\ldots,m_k,n_k)$ and $(m'_1,n'_1,m'_2,n'_2,\ldots,\allowbreak m'_{k+1},\allowbreak n'_{k+1})$ are
  $(m,n)$-\textit{pro\-jec\-tive factorization sequences}.

  Furthermore, we let $\tau(m,n)$ denote the \textit{number of distinct $(m,n)$-pro\-jec\-tive factorization sequences}.
\end{definition}

The next result provides the number of distinct $C_4$-face-magic projective labelings on $\mathcal{P}_{m,n}$ having
$C_4$-face-magic value $2mn+2$ for distinct odd integers $m$ and $n$.

\begin{theorem}[\cite{Curran2}, Theorem 47] \label{thmNumberOfCentrallyBalancedC4FaceMagicLabelingsOnRectangle}
  Let $m \geqslant 3$ and $n \geqslant 3$ be distinct odd integers.
  Then the number of distinct $C_4$-face-magic projective labelings on $\mathcal{P}_{m,n}$  having $C_4$-face-magic value $2mn+2$
  (up to symmetries on the projective plane) is
  \begin{equation*}
    \bigl( \tau(m,n) + \tau(n,m) \bigr) \, 2^{m/2 +n/2-3} \, (\tfrac{m-1}{2})! (\tfrac{n-1}{2})!.
  \end{equation*}
\end{theorem}

The number of distinct $C_4$-face-magic projective labelings on $\mathcal{P}_{m,m}$ having
$C_4$-face-magic value $2m^2+2$ for an odd integer $m$ is stated below.

\begin{theorem}[\cite{Curran2}, Theorem 48] \label{thmNumberOfCentrallyBalancedC4FaceMagicLabelingsOnSquare}
  Let $m \geqslant 3$ be an odd integer.
  Then the number of distinct $C_4$-face-magic projective labelings on $\mathcal{P}_{m,m}$ having $C_4$-face-magic value $2m^2+2$
  (up to symmetries on the projective plane) is
  \begin{equation*}
    \tau(m,m) \, 2^{m-3} \, \bigl( (\tfrac{m-1}{2})! \bigr)^2.
  \end{equation*}
\end{theorem}

We now direct our attention to $C_4$-face-magic labelings on $\mathcal{P}_{m,n}$  with $C_4$-face-magic value $2mn+1$ or $2mn+3$.

\begin{definition} \label{defnOrderComplementLabeling}
  Let $X=\{ x_{i,j} : (i,j)\in V(\mathcal{P}_{m,n})\}$ be a $C_4$-face-magic projective labeling on $\mathcal{P}_{m,n}$.
  We define a labeling $Y$ on $\mathcal{P}_{m,n}$ given by
  \begin{equation*}
   y_{i,j} = mn+1 -x_{i,j} \mbox{ \ \ \ for all } (i,j)\in V(\mathcal{P}_{m,n}).
  \end{equation*}
  We say that $Y$ is the \textit{order plus one complement of $X$ labeling} on $\mathcal{P}_{m,n}$.
  We call the transformation $\mathcal{C}(X)=Y$ the \textit{order plus one complement labeling transformation}.
\end{definition}

\begin{proposition} \label{propOrderComplemtLabeling}
  Let $X=\{ x_{i,j} : (i,j)\in V(\mathcal{P}_{m,n})\}$ be a $C_4$-face-magic projective labeling on $\mathcal{P}_{m,n}$
  with $C_4$-face-magic value $S$.
  Then the order plus one complement of $X$ labeling $Y=\mathcal{C}(X)$ on $\mathcal{P}_{m,n}$ is a $C_4$-face-magic projective labeling on $\mathcal{P}_{m,n}$
  with $C_4$-face-magic value $4mn+4-S$.
\end{proposition}

\begin{proof}
  Let $(i_1,j_1)$, $(i_2,j_2)$, $(i_3,j_3)$, and $(i_4,j_4)$ be the vertices of any $C_4$-face on $\mathcal{P}_{m,n}$. Then
  \begin{equation*}
    \sum_{k=1}^4 x_{i_k,j_k} = S.
  \end{equation*}
  Thus, the order complement of $X$ labeling $Y=\mathcal{C}(X)$ satisfies
   \begin{equation*}
    \sum_{k=1}^4 y_{i_k,j_k} =   \sum_{k=1}^4 (mn+1 -x_{i_k,j_k}) = 4mn+4-S.
  \end{equation*}
  Since $h:\{1,2,\ldots,mn\}\rightarrow \{1,2,\ldots,mn\}$ defined by $h(x)=mn+1-x$ is a bijection,
  so is $f:V(P_m \times P_n) \rightarrow \{1,2,\ldots,mn\}$ defined by $f(i,j)= h( x_{i,j})$.
\end{proof}

\begin{remark} \label{remOrderComplementLabelingCorrespondence}
  We observe that $\mathcal{C}(\mathcal{C}(X))=X$ for all $C_4$-face-magic projective labelings
  $X$ on $\mathcal{P}_{m,n}$.
  Thus, $\mathcal{C}$ is a one-to-one correspondence between
  $C_4$-face-magic projective labelings on $\mathcal{P}_{m,n}$ with $C_4$-face-magic value $S$ and
  $C_4$-face-magic projective labelings on $\mathcal{P}_{m,n}$ with $C_4$-face-magic value $4mn+4-S$.

  Hence, for odd integers $m\geqslant 3$ and $n \geqslant 3$,
  $\mathcal{C}$ is a one-to-one correspondence between
  $C_4$-face-magic projective labelings on $\mathcal{P}_{m,n}$ with $C_4$-face-magic value $2mn+1$ and
  those with $C_4$-face-magic value $2mn+3$.
\end{remark}

\section{Bicentrally balanced $C_4$-face-magic labelings on $\mathcal{P}_{m,n}$} \label{S:NonstandardSumProjectiveGrid}

\begin{notation}
Throughout this section, we assume that both $m\geqslant 3$ and $n\geqslant 3$ are odd integers.
We write $m=2m_0 +1$ and $n=2n_0 +1$ for positive integers $m_0$ and $n_0$.
For any positive integer $N$, we let $N^{+}=N+1$.
In particular, we have $m_0^{+}=m_0 +1$ and $n_0^{+}=n_0 +1$.
\end{notation}

\begin{notation}
  We refer to the vertex $(\tfrac{1}{2}(m+1),\tfrac{1}{2}(n+1))=(m_0^{+},n_0^{+})$ as the \textit{center} of the projective grid graph $\mathcal{P}_{m,n}$.
  The graph automorphisms of $\mathcal{P}_{m,n}$ that are induced by a homeomorphism of the projective plane are
  described in relation to the center of $\mathcal{P}_{m,n}$.
  We let $R_{\theta}$ denote the rotation by $\theta$ degrees in the counter-clockwise direction about the center.
  The symmetry $H$ ($V$) is the reflection about the center row (column).
  Since the corner vertices of $\mathcal{P}_{m,n}$ are the only vertices incident to a double edge,
  a symmetry of $\mathcal{P}_{m,n}$ sends each corner vertex to another corner vertex.
  Thus, for distinct integers $m$ and $n$, the set of symmetries on $\mathcal{P}_{m,n}$ is $\{ R_{0}, R_{180}, H, V\}$.
  We let $D_{+}$ ($D_{-}$) denote the reflection about the diagonal with positive (negative) slope passing through the center.
  When $m=n$, the set of symmetries on $\mathcal{P}_{m,m}$ is
  $D_4=\{R_{0}, R_{90}, R_{180}, R_{270}, H, V, D_{+}, D_{-}\}$.
\end{notation}

\begin{definition} \label{defnBicentrallyBalanced}
Let $X=\{x_{i,j}: (i,j) \in V(\mathcal{P}_{m,n}) \}$ be a $C_4$-face-magic projective labeling on $\mathcal{P}_{m,n}$ with
$C_4$-face value $S=2mn+3$.
For all $(i,j) \in V(\mathcal{P}_{m,n})$, let
\begin{align*}
  S(i,j) &= \tfrac{1}{2}mn+ \tfrac{3}{2} \mbox{ if } i+j \mbox{ is even, and } \\
   S(i,j) &= \tfrac{3}{2}mn+ \tfrac{3}{2} \mbox{ if } i+j \mbox{ is odd.}
\end{align*}
We say that $X$ is \emph{bicentrally balanced} if, for all $(i,j) \in V(\mathcal{P}_{m,n})$,
\begin{equation*}
  x_{i,j} + x_{m+1-i,n+1-j} = S(i,j).
\end{equation*}
\end{definition}

\begin{remark} \label{remBicentrallyBalancedHasValue2MNPlus3}
Suppose $X=\{x_{i,j}: (i,j) \in V(\mathcal{P}_{m,n}) \}$ is a $C_4$-face-magic labeling on $\mathcal{P}_{m,n}$ that is bicentrally balanced.
Then the $C_4$-face-magic value $S$ of $X$ is
\begin{equation*}
  S = x_{1,1} + x_{m,n} + x_{2,1} + x_{m-1,n}=  S(1,1) + S(2,1) = 2mn+3.
\end{equation*}
\end{remark}

We observe that the $C_4$-face-magic projective labeling on $\mathcal{P}_{5,5}$ in Fig. \ref{figFiveByFiveProjectiveGridGraphLabeling} is bicentrally balanced.

\begin{lemma} \label{lemmaProjectiveGridSemiC4FaceValueBalance}
Suppose $m \geqslant 3$ and  $n \geqslant 3$ are odd integers.
Let $X = \{ x_{i,j} : (i,j)\in V(\mathcal{P}_{m,n}) \}$ be a $C_4$-face-magic projective labeling on $\mathcal{P}_{m,n}$
with $C_4$-face-magic value $S=2mn+3$.
Then $X$ is bicentrally balanced. Furthermore, we have $x_{m_0^{+},n_0^{+}}= \tfrac{1}{2}S(m_0^{+},n_0^{+})$.
Thus, $x_{m_0^{+},n_0^{+}}=\tfrac{1}{4}mn+ \tfrac{3}{4}$ if $m^{+}_0 + n^{+}_0$ is even, or
$x_{m_0^{+},n_0^{+}}=\tfrac{3}{4}mn+ \tfrac{3}{4}$ if $m^{+}_0 + n^{+}_0$ is odd.
\end{lemma}

\begin{proof}
By Lemma \ref{lemmaOddByOddProjectiveGridC4FaceValue}, the digons formed by the
vertex sets $\{ (1,1), (m,n)\}$ and $\{(m,1),(1,n)\}$ have face values
\begin{align*}
  D_1 &= x_{1,1} + x_{m,n} = \tfrac{1}{2}mn+ \tfrac{3}{2} \mbox{ and } \\
  D_2 &=  x_{m,1} + x_{1,n} = \tfrac{1}{2}mn+ \tfrac{3}{2}, \mbox{ respectively.}
\end{align*}

Suppose that for some integer $1 \leqslant i < m$,
\begin{equation*}
  x_{i,1} + x_{m+1-i,n} =S(i,1).
\end{equation*}
See Definition \ref{defnBicentrallyBalanced} for the definition of $S(i,j)$.
Since
\begin{equation*}
  x_{i,1} + x_{i+1,1} + x_{m+1-i,n}  + x_{m-i,n} =S,
\end{equation*}
we have
\begin{equation*}
  x_{i+1,1} + x_{m-i,n} =S - S(i,1)= S(i+1,1).
\end{equation*}
Hence,
\begin{equation*}
  x_{i,1} + x_{m+1-i,n} =S(i,1)
\end{equation*}
 for all $1 \leqslant i \leqslant m$.
A similar argument shows that
\begin{equation*}
  x_{1,j} + x_{m,n+1-j} =S(1,j)
\end{equation*}
for all $1 \leqslant j \leqslant n$.

We use induction to show that $x_{i,j} + x_{m+1-i,n+1-j} = S(i,j)$ for all $(i,j)\in V\mathcal{P}_{m,n})$.
Suppose there exist integers $1 < i <m$ and $1 < j < n$ such that
\begin{enumerate}
  \item for all $1\leqslant i' <i$ and $1\leqslant j' \leqslant n$, $x_{i',j'} + x_{m+1-i',n+1-j'} = S(i',j')$, and
  \item for all $1\leqslant j' < j$, $x_{i,j'} + x_{m+1-i,n+1-j'} = S(i,j')$.
\end{enumerate}
We need to show that $x_{i,j} + x_{m+1-i,n+1-j} = S(i,j)$.
When we add the two $C_4$-face-values
\begin{align*}
  x_{i-1,j-1} &+ x_{i-1,j} + x_{i,j-1} + x_{i,j} = S \\
  &\mbox{ and } \\
  x_{m+2-i,n+2-j} &+ x_{m+2-i,n+1-j} + x_{m+1-i,n+2-j} + x_{m+1-i,n+1-j} = S,
\end{align*}
we obtain
\begin{align*}
  (x_{i-1,j-1} + x_{m+2-i,n+2-j}) &+ (x_{i-1,j}  + x_{m+2-i,n+1-j} ) \\
  + (x_{i,j-1}+ x_{m+1-i,n+2-j}) &+ (x_{i,j}+ x_{m+1-i,n+1-j}) = 2S.
\end{align*}
Since
\begin{align*}
   x_{i-1,j-1} + x_{m+2-i,n+2-j} &= S(i-1,j-1), \\
   x_{i-1,j}  + x_{m+2-i,n+1-j} &= S(i-1,j), \mbox{ and} \\
   x_{i,j-1}+ x_{m+1-i,n+2-j}  &= S(i,j-1),
\end{align*}
we have
\begin{equation*}
  S(i-1,j-1) + S(i-1,j) + S(i,j-1) + (x_{i,j}+ x_{m+1-i,n+1-j}) = 2S.
\end{equation*}
Thus
\begin{equation*}
x_{i,j} + x_{m+1-i,n+1-j} =S(i,j).
\end{equation*}

Since
\begin{equation*}
  2x_{m_0^{+},n_0^{+}}= x_{m_0^{+},n_0^{+}} + x_{m+1-m_0^{+},n+1-n_0^{+}} = S(m^{+}_0,n^{+}_0),
\end{equation*}
we have
\begin{equation*}
  x_{m_0^{+},n_0^{+}}= \tfrac{1}{2}S(m^{+}_0,n^{+}_0).
\end{equation*}
Thus, $x_{m_0^{+},n_0^{+}}=\tfrac{1}{4}mn+ \tfrac{3}{4}$ if $m^{+}_0 + n^{+}_0$ is even, or
$x_{m_0^{+},n_0^{+}}=\tfrac{3}{4}mn+ \tfrac{3}{4}$ if $m^{+}_0 + n^{+}_0$ is odd.
\end{proof}

\subsection{Structure of a bicentrally balanced $C_4$-face-magic labeling}

\begin{lemma} \label{lemmaStructureOfBicentrallyBalancedLabeling}
   Let $X=\{x_{i,j}:(i,j)\in V(\mathcal{P}_{m,n})\}$ be a bicentrally balanced $C_4$-face-magic projective labeling on $\mathcal{P}_{m,n}$.
   For $1 \leqslant j \leqslant n_0$, let
   \begin{equation*}
     a_j = x_{1,j} +x_{1,j+1}.
   \end{equation*}
   Then,
   \begin{enumerate}
     \item for all $1 \leqslant i \leqslant m_0$ where $i$ is odd, and $1 \leqslant j \leqslant n_0$, we have
     \begin{align*}
       x_{i,j} + x_{i,j+1} &= a_j,                      &x_{i,n+1-j} + x_{i,n-j} &= S- a_j,  \\
       x_{m+1-i,j} + x_{m+1-i,j+1} &= a_j,  \mbox{and } &x_{m+1-i,n+1-j} + x_{m+1-i,n-j} &= S- a_j, \mbox{\ \ and }
     \end{align*}
     \item for all $1 \leqslant i \leqslant m_0$ where $i$ is even, and $1 \leqslant j \leqslant n_0$, we have
     \begin{align*}
       x_{i,j} + x_{i,j+1} &= S- a_j,                      &x_{i,n+1-j} + x_{i,n-j} &= a_j,  \\
       x_{m+1-i,j} + x_{m+1-i,j+1} &= S-a_j,  \mbox{and } &x_{m+1-i,n+1-j} + x_{m+1-i,n-j} &= a_j.
     \end{align*}
   \end{enumerate}
\end{lemma}

\begin{proof}
  When we equate the two $C_4$-face sums
  \begin{align*}
    x_{i,j} + x_{i,j+1} +  x_{i+1,j} + x_{i+1,j+1}  &= S \mbox{ \ and} \\
     x_{i+1,j} + x_{i+1,j+1} + x_{i+2,j} + x_{i+2,j+1} &=S,
  \end{align*}
  we obtain
  \begin{equation} \label{eqnEveryOtherPair}
     x_{i,j} + x_{i,j+1}  =  x_{i+2,j} + x_{i+2,j+1}.
  \end{equation}
  By (\ref{eqnEveryOtherPair}), for all $1 \leqslant i \leqslant m_0$ where $i$ is odd, and $1\leqslant j \leqslant n_0$, we have
  \begin{equation*}
    x_{i,j} + x_{i,j+1} = a_j \mbox{ \ and \ }  x_{m+1-i,j} + x_{m+1-i,j+1} = a_j.
  \end{equation*}
  Since
  \begin{equation*}
    a_j + x_{2,j} + x_{2,j+1} = x_{1,j} +x_{1,j+1} + x_{2,j} +x_{2,j+1} = S,
  \end{equation*}
  we have
  \begin{equation*}
    x_{2,j} + x_{2,j+1} = S- a_j
  \end{equation*}
  for all $1 \leqslant j \leqslant n_0$.
  Also by (\ref{eqnEveryOtherPair}), for all $1 \leqslant i \leqslant m_0$ where $i$ is even, and $1\leqslant j \leqslant n_0$, we have
  \begin{equation*}
    x_{i,j} + x_{i,j+1} = S-a_j \mbox{ \ and \ }  x_{m+1-i,j} + x_{m+1-i,j+1} = S-a_j.
  \end{equation*}

  Since
  \begin{equation*}
    a_j + x_{1,n+1-j} + x_{1,n-j} = x_{m,j} +x_{m,j+1} + x_{1,n+1-j} +x_{1,n-j} = S,
  \end{equation*}
  we have
  \begin{equation*}
    x_{1,n+1-j} +x_{1,n-j} = S- a_j
  \end{equation*}
  for all $1 \leqslant j \leqslant n_0$.
  Thus, by (\ref{eqnEveryOtherPair}), for all $1 \leqslant i \leqslant m_0$ where $i$ is odd, and $1\leqslant j \leqslant n_0$, we have
  \begin{equation*}
    x_{i,n+1-j} + x_{i,n-j} = S - a_j \mbox{ \ and \ }  x_{m+1-i,n+1-j} + x_{m+1-i,n-j} = S - a_j.
  \end{equation*}
  Since
  \begin{equation*}
    (S -a_j) + x_{2,n+1-j} + x_{2,n-j} =  x_{1,n+1-j} +x_{1,n-j} + x_{2,n+1-j} + x_{2,n-j}  = S,
  \end{equation*}
  we have
  \begin{equation*}
    x_{2,n+1-j} + x_{2,n-j}  = a_j
  \end{equation*}
  for all $1 \leqslant j \leqslant n_0$.
  Hence, by (\ref{eqnEveryOtherPair}), for all $1 \leqslant i \leqslant m_0$ where $i$ is even, and $1\leqslant j \leqslant n_0$, we have
  \begin{equation*}
    x_{i,n+1-j} + x_{i,n-j} = a_j \mbox{ \ and \ }  x_{m+1-i,n+1-j} + x_{m+1-i,n-j} = a_j.
  \end{equation*}
\end{proof}

\subsection{Row and column permutations on a bicentrally balanced labeling}

\begin{definition} \label{defnElementaryProjectivePermutationColumns}
  Let $X=\{x_{i,j}:(i,j)\in V(\mathcal{P}_{m,n})\}$ be a bicentrally balanced $C_4$-face-magic projective labeling on $\mathcal{P}_{m,n}$.
  Let $\eta$ be a permutation on the set $\{1,2,\ldots,m_0\}$ such that $\eta(i)\equiv i \pmod{2}$ for all $1 \leqslant i \leqslant m_0$.
  We define a labeling on $\mathcal{P}_{m,n}$, $Z=\{z_{i,j}:(i,j)\in V(\mathcal{P}_{m,n})\}$,
  such that for all $1 \leqslant i \leqslant m_0$ and $1 \leqslant j \leqslant n$, we have
    \begin{align*}
      z_{i,j} &= x_{\eta(i),j},  \\
      z_{m_0^{+},j} &= x_{m_0^{+},j}, \mbox{ and}\\
       z_{m+1-i,j} &= x_{m+1-\eta(i),j}.
    \end{align*}
    We let $\mathcal{E}_{\eta}$ denote the labeling operation given by $\mathcal{E}_{\eta}(X)=Z$.
\end{definition}

\begin{lemma}  \label{lemmaElementaryProjectivePermutationColumns}
  Let $X=\{x_{i,j}:(i,j)\in V(\mathcal{P}_{m,n})\}$ be a bicentrally balanced $C_4$-face-magic projective labeling on $\mathcal{P}_{m,n}$, and
  let $\eta$ be a permutation on the set $\{1,2,\ldots,m_0\}$ such that $\eta(i)\equiv i \pmod{2}$ for all $1 \leqslant i \leqslant m_0$.
  Let $\mathcal{E}_{\eta}$ be the labeling operation defined in Definition \ref{defnElementaryProjectivePermutationColumns}.
  Then the labeling $Z=\mathcal{E}_{\eta}(X)$ is a bicentrally balanced $C_4$-face-magic projective labeling on $\mathcal{P}_{m,n}$.
\end{lemma}

\begin{proof}
 We first verify that $Z$ is bicentrally balanced.
 Suppose that $1 \leqslant i \leqslant m_0$ and $1 \leqslant j \leqslant n$.
 Since $\eta(i)-i$ is even, we have
 \begin{equation*}
   z_{i,j} + z_{m+1-i,n+1-j} = x_{\eta(i),j} + x_{m+1-\eta(i),n+1-j}  =S(i,j).
 \end{equation*}
 Furthermore, we have
  \begin{equation*}
   z_{m_0^{+},j} + z_{m+1-m_0^{+},n+1-j} = x_{m_0^{+},j} + x_{m+1-m_0^{+},n+1-j}  = S(m^{+}_0,j).
 \end{equation*}

  Next, we show that $Z$ is a $C_4$-face-magic projective labeling on $\mathcal{P}_{m,n}$.
  For all $1 \leqslant i <m$ and $1 \leqslant j <n$, by Lemma \ref{lemmaStructureOfBicentrallyBalancedLabeling}, we have
  \begin{equation*}
    z_{i,j} + z_{i,j+1} = x_{i,j} + x_{i,j+1} \mbox{ \ and \ }  z_{i+1,j} + z_{i+1,j+1} = x_{i+1,j} + x_{i+1,j+1}.
  \end{equation*}
  Thus
  \begin{align*}
    z_{i,j} &+ z_{i,j+1} + z_{i+1,j} + z_{i+1,j+1} \\
    & = x_{i,j} + x_{i,j+1} + x_{i+1,j} + x_{i+1,j+1} =S.
  \end{align*}
  Since $Z$ is bicentrally balanced, for $1 \leqslant i <m$, we have
  \begin{equation*}
    z_{i,n} +z_{m+1-i,1} + z_{i+1,n} + z_{m-i,1} = S(i,n) + S(i+1,n) = S.
  \end{equation*}
  Also, since $Z$ is bicentrally balanced, for $1 \leqslant j < n$, we have
  \begin{equation*}
    z_{m,j} +z_{1,n+1-j} + z_{m,j+1} + z_{1,n-j} = S(m,j) + S(m,j+1) = S.
  \end{equation*}
\end{proof}

\begin{definition} \label{defnElementaryProjectivePermutationRows}
  Let $X=\{x_{i,j}:(i,j)\in V(\mathcal{P}_{m,n})\}$ be a bicentrally balanced $C_4$-face-magic projective labeling on $\mathcal{P}_{m,n}$.
  Let $\kappa$ be a permutation on the set $\{1,2,\ldots,n_0\}$ such that $\kappa(j)\equiv j \pmod{2}$ for all $1 \leqslant j \leqslant n_0$.
  We define a labeling on $\mathcal{P}_{m,n}$, $Z=\{z_{i,j}:(i,j)\in V(\mathcal{P}_{m,n})\}$,
  such that for all $1 \leqslant i \leqslant m$ and $1 \leqslant j \leqslant n_0$, we have
    \begin{align*}
      z_{i,j} &= x_{i,\kappa(j)},  \\
      z_{i,n_0^{+}} &= x_{i,n_0^{+}}, \mbox{ and } \\
      z_{i,n+1-j} &= x_{i,n+1-\kappa(j)}.
    \end{align*}
    We let $\mathcal{E}_{\kappa}$ denote the labeling operation given by $\mathcal{E}_{\kappa}(X)=Z$.
\end{definition}

\begin{lemma}  \label{lemmaElementaryProjectivePermutationRows}
  Let $X=\{x_{i,j}:(i,j)\in V(\mathcal{P}_{m,n})\}$ be a bicentrally balanced $C_4$-face-magic projective labeling on $\mathcal{P}_{m,n}$, and
  let $\kappa$ be a permutation on the set $\{1,2,\ldots,n_0\}$ such that $\kappa(j)\equiv j \pmod{2}$ for all $1 \leqslant j \leqslant n_0$.
  Let $\mathcal{E}_{\kappa}$ be the labeling operation defined in Definition \ref{defnElementaryProjectivePermutationRows}.
  Then the labeling $Z=\mathcal{E}_{\kappa}(X)$ is a bicentrally balanced $C_4$-face-magic projective labeling on $\mathcal{P}_{m,n}$.
\end{lemma}

The proof of Lemma \ref{lemmaElementaryProjectivePermutationRows} is similar to the proof of Lemma \ref{lemmaElementaryProjectivePermutationColumns}.

\begin{definition} \label{defnElementaryProjectiveSwapColumns}
  Let $X=\{x_{i,j}:(i,j)\in V(\mathcal{P}_{m,n})\}$ be a bicentrally balanced $C_4$-face-magic projective labeling on $\mathcal{P}_{m,n}$.
  Let $\alpha: \{1,2,\ldots,m_0\}\rightarrow \{0,1\}$.
  We define a labeling on $\mathcal{P}_{m,n}$, $Z=\{z_{i,j}:(i,j)\in V(\mathcal{P}_{m,n})\}$,
  such that for all $1 \leqslant i \leqslant m_0$ and $1 \leqslant j \leqslant n$, we have
    \begin{align*}
      z_{i,j} &= x_{(1-\alpha(i))i +\alpha(i)(m+1-i),j}, \mbox{ \ \ and} \\
      z_{m+1-i,j} &= x_{\alpha(i)i + (1- \alpha(i)) (m+1-i),j}.
    \end{align*}
    We let $\mathcal{E}_{\alpha}$ denote the labeling operation given by $\mathcal{E}_{\alpha}(X)=Z$.
    The labeling operation $\mathcal{E}_{\alpha}$ has the effect of keeping the labelings on the vertices of columns $i$ and $m+1-i$ the same if $\alpha(i)=0$ and
    swapping the labelings on the vertices of column $i$ with those of column $m+1-i$ if $\alpha(i)=1$.
\end{definition}

\begin{lemma}  \label{lemmaElementaryProjectiveSwapColumns}
  Let $X=\{x_{i,j}:(i,j)\in V(\mathcal{P}_{m,n})\}$ be a bicentrally balanced $C_4$-face-magic projective labeling on $\mathcal{P}_{m,n}$, and
  let $\alpha: \{1,2,\ldots,m_0\}\rightarrow \{0,1\}$.
  Let $\mathcal{E}_{\alpha}$ be the labeling operation defined in Definition \ref{defnElementaryProjectiveSwapColumns}.
  Then the labeling $Z=\mathcal{E}_{\alpha}(X)$ is a bicentrally balanced $C_4$-face-magic projective labeling on $\mathcal{P}_{m,n}$.
\end{lemma}

\begin{proof}
  First, we show that $Z$ is bicentrally balanced.
  Suppose $\alpha(i)=0$. Then
  \begin{equation*}
    z_{i,j} = x_{i,j} \mbox{ \ \ and \ \ } z_{m+1-i,j} = x_{m+1-i,j}.
  \end{equation*}
  Thus
  \begin{equation*}
    z_{i,j} +  z_{m+1-i,n+1-j} =  x_{i,j} +  x_{m+1-i,n+1-j} = S(i,j).
  \end{equation*}
  Suppose $\alpha(i)=1$. Then
  \begin{equation*}
    z_{i,j} = x_{m+1-i,j} \mbox{ \ \ and \ \ } z_{m+1-i,j} = x_{i,j}.
  \end{equation*}
  Thus
  \begin{equation*}
    z_{i,j} +  z_{m+1-i,n+1-j} =  x_{m+1-i,j} +  x_{i,n+1-j} = S(m+1-i,j) =S(i,j)
  \end{equation*}
  since $m+1-i \equiv i \pmod{2}$.
  The proof that $Z$ is a $C_4$-face-magic projective labeling on $\mathcal{P}_{m,n}$ with $C_4$-face-magic value $2mn+3$ is similar to
  that in the proof of Lemma \ref{lemmaElementaryProjectivePermutationColumns}.
\end{proof}

\begin{definition} \label{defnElementaryProjectiveSwapRows}
  Let $X=\{x_{i,j}:(i,j)\in V(\mathcal{P}_{m,n})\}$ be a bicentrally balanced $C_4$-face-magic projective labeling on $\mathcal{P}_{m,n}$.
  Let $\delta: \{1,2,\ldots,n_0\}\rightarrow \{0,1\}$.
  We define a labeling on $\mathcal{P}_{m,n}$, $Z=\{z_{i,j}:(i,j)\in V(\mathcal{P}_{m,n})\}$,
  such that for all $1 \leqslant i \leqslant m$ and $1 \leqslant j \leqslant n_0$, we have
    \begin{align*}
      z_{i,j} &= x_{i, (1-\delta(j))j +\delta(j)(n+1-j)}, \mbox{ \ \ and} \\
      z_{i,n+1-j} &= x_{i,\delta(j)j + (1- \delta(j)) (n+1-j)}.
    \end{align*}
    We let $\mathcal{E}_{\delta}$ denote the labeling operation given by $\mathcal{E}_{\delta}(X)=Z$.
    The labeling operation $\mathcal{E}_{\delta}$ has the effect of keeping the labelings on the vertices of rows $j$ and $n+1-j$ the same if $\delta(j)=0$ and
    swapping the labelings on the vertices of row $j$ with those of row $n+1-j$ if $\delta(j)=1$.
\end{definition}

\begin{lemma}  \label{lemmaElementaryProjectiveSwapRows}
  Let $X=\{x_{i,j}:(i,j)\in V(\mathcal{P}_{m,n})\}$ be a bicentrally balanced $C_4$-face-magic projective labeling on $\mathcal{P}_{m,n}$, and
  let $\delta: \{1,2,\ldots,n_0\}\rightarrow \{0,1\}$.
  Let $\mathcal{E}_{\delta}$ be the labeling operation defined in Definition \ref{defnElementaryProjectiveSwapRows}.
  Then the labeling $Z=\mathcal{E}_{\delta}(X)$ is a bicentrally balanced $C_4$-face-magic projective labeling on $\mathcal{P}_{m,n}$.
\end{lemma}

The proof of Lemma \ref{lemmaElementaryProjectiveSwapRows} is similar to the proof of Lemma \ref{lemmaElementaryProjectiveSwapColumns}.

\begin{definition}
  We call each of the labeling operations $\mathcal{E}_{\eta}$ in Definition \ref{defnElementaryProjectivePermutationColumns},
  $\mathcal{E}_{\kappa}$ in Definition \ref{defnElementaryProjectivePermutationRows},
  $\mathcal{E}_{\alpha}$ in Definition \ref{defnElementaryProjectiveSwapColumns}, and
  $\mathcal{E}_{\delta}$ in Definition \ref{defnElementaryProjectiveSwapRows} an
  \textit{elementary projective labeling operation}.
\end{definition}

\begin{definition}
  We say that two bicentrally balanced $C_4$-face-magic labelings on $\mathcal{P}_{m,n}$ are
  \textit{projective labeling equivalent} if one labeling can be obtained from the other by applying a sequence
  of elementary projective labeling operations to it.
\end{definition}

\subsection{Standard bicentrally balanced labeling}

Given a bicentrally balanced $C_4$-face-magic projective labeling $X$ on $\mathcal{P}_{m,n}$, the next theorem identifies
a canonical bicentrally balanced $C_4$-face-magic projective labeling on $\mathcal{P}_{m,n}$ that is
projective labeling equivalent to $X$.

\begin{theorem} \label{thmStandardProjectiveLabeling}
  Let $X=\{x_{i,j}:(i,j)\in V(\mathcal{P}_{m,n})\}$ be a bicentrally balanced $C_4$-face-magic projective labeling on $\mathcal{P}_{m,n}$.
  Then there is a unique bicentrally balanced $C_4$-face-magic projective labeling $Z=\{z_{i,j}:(i,j)\in V(\mathcal{P}_{m,n})\}$
  on $\mathcal{P}_{m,n}$ that is projective labeling equivalent to $X$ such that
  \begin{itemize}
    \item $z_{i,n_0^{+}} <  z_{i+2,n_0^{+}}$ for all $1 \leqslant i \leqslant m-2$ and $i+ n^{+}_0$ is even,
    \item $z_{i,n_0^{+}} >  z_{i+2,n_0^{+}}$ for all $1 \leqslant i \leqslant m-2$ and $i+ n^{+}_0$ is odd,
    \item $z_{m_0^{+},j} <  z_{m_0^{+},j+2}$ for all $1 \leqslant j \leqslant n-2$ and $m^{+}_0 +j$ is even, and
    \item $z_{m_0^{+},j} >  z_{m_0^{+},j+2}$ for all $1 \leqslant j \leqslant n-2$ and $m^{+}_0 +j$ is odd.
  \end{itemize}
\end{theorem}

\begin{proof}
  By Lemma \ref{lemmaProjectiveGridSemiC4FaceValueBalance}, we have $x_{m_0^{+},n_0^{+}}=\tfrac{1}{2} S(m^{+}_0,n^{+}_0)$.
  It is easy to check that this value remains the same regardless of the elementary projective labeling operation that we apply to $X$.
  Since $X$ is bicentrally balanced, for all $1 \leqslant i \leqslant m_0$, we have
  \begin{equation*}
    x_{i,n_0^{+}} + x_{m+1-i,n_0^{+}} =S(i,n^{+}_0).
  \end{equation*}
  Thus, either  $x_{i,n_0^{+}} < \tfrac{1}{2}S(i,n^{+}_0)$ or $x_{m+1-i,n_0^{+}} < \tfrac{1}{2}S(i,n^{+}_0)$.
  We define a function $\alpha:\{1,2,\allowbreak \ldots,\allowbreak m_0\}\allowbreak \rightarrow \allowbreak \{0,1\}$ as follows.
  For each $1 \leqslant i \leqslant m_0$, let
  \begin{equation*}
    \alpha(i) = \left\{\begin{array}{ll}
              0, &\mbox{if \ } i+n^{+}_0 \mbox{ is even and } x_{i,n_0^{+}} < \tfrac{1}{2}S(i,n^{+}_0), \\
              1, &\mbox{if \ } i+n^{+}_0 \mbox{ is even and } x_{i,n_0^{+}} > \tfrac{1}{2}S(i,n^{+}_0), \\
              0, &\mbox{if \ } i+n^{+}_0 \mbox{ is odd and } x_{i,n_0^{+}} > \tfrac{1}{2}S(i,n^{+}_0), \mbox{ and } \\
              1, &\mbox{if \ } i+n^{+}_0 \mbox{ is odd and } x_{i,n_0^{+}} < \tfrac{1}{2}S(i,n^{+}_0).
\end{array} \right.
  \end{equation*}
  By Lemma \ref{lemmaElementaryProjectiveSwapColumns},  $\mathcal{E}_{\alpha}(X)$ is a bicentrally balanced $C_4$-face-magic projective labeling on $\mathcal{P}_{m,n}$.
  Replace $X$ with $\mathcal{E}_{\alpha}(X)$. Then $X$ satisfies,
  for all $1 \leqslant i \leqslant m_0$,
  \begin{align*}
    x_{i,n_0^{+}} &<\tfrac{1}{2}S(i,n^{+}_0) \mbox{ and } x_{m+1-i,n_0^{+}} > \tfrac{1}{2}S(i,n^{+}_0) \mbox{ if } i+n^{+}_0 \mbox{ is even, and } \\
     x_{i,n_0^{+}} &>\tfrac{1}{2}S(i,n^{+}_0) \mbox{ and } x_{m+1-i,n_0^{+}} < \tfrac{1}{2}S(i,n^{+}_0) \mbox{ if } i+n^{+}_0 \mbox{ is odd.}
  \end{align*}
  Choose a permutation $\eta$ of $\{1,2,\ldots,m_0\}$ with $\eta(i)\equiv i \pmod{2}$ for all $1 \leqslant i \leqslant m_0$
  such that for all  $1 \leqslant i \leqslant m_0 -2$,
  \begin{align*}
    x_{\eta(i),n_0^{+}} &< x_{\eta(i+2),n_0^{+}} \mbox{ \ if }  i+n^{+}_0 \mbox{ is even, and } \\
    x_{\eta(i),n_0^{+}} &> x_{\eta(i+2),n_0^{+}} \mbox{ \ if }  i+n^{+}_0 \mbox{ is odd.}
  \end{align*}
  By Lemma \ref{lemmaElementaryProjectivePermutationColumns},
  $\mathcal{E}_{\eta}(X)$ is a  bicentrally balanced $C_4$-face-magic projective labeling on $\mathcal{P}_{m,n}$.
  Replace $X$ with $\mathcal{E}_{\eta}(X)$. Then $X$ satisfies,
  for all $1 \leqslant i \leqslant m-2$,
  \begin{align*}
    x_{i,n^{+}_0}   &<  x_{i+2,n^{+}_0}  \mbox{ \ if }  i+n^{+}_0 \mbox{ is even, and } \\
     x_{i,n^{+}_0}   &>  x_{i+2,n^{+}_0}  \mbox{ \ if }  i+n^{+}_0 \mbox{ is odd.}
  \end{align*}

  A similar argument allows us to choose a function $\delta:\{1,2,\ldots,n_0\}\rightarrow \{0,1\}$ and a permutation $\kappa$ on $\{1,2,\ldots,n_0\}$
  with $\kappa(j)\equiv j \pmod{2}$ for all $1 \leqslant j \leqslant n_0$
  such that $Z=\mathcal{E}_{\kappa}( \mathcal{E}_{\delta}(X))$
  is a  bicentrally balanced $C_4$-face-magic projective labeling on $\mathcal{P}_{m,n}$ and
  \begin{align*}
    z_{i,n_0^{+}} &<  z_{i+2,n_0^{+}} \mbox{ \ for all } 1 \leqslant i \leqslant m-2 \mbox{ and } i+ n^{+}_0 \mbox{ is even,}  \\
    z_{i,n_0^{+}} &>  z_{i+2,n_0^{+}} \mbox{ \ for all } 1 \leqslant i \leqslant m-2 \mbox{ and } i+ n^{+}_0 \mbox{ is odd,}  \\
    z_{m_0^{+},j} &<  z_{m_0^{+},j+2} \mbox{ \ for all } 1 \leqslant j \leqslant n-2 \mbox{ and } m^{+}_0 +j \mbox{ is even, and} \\
    z_{m_0^{+},j} &>  z_{m_0^{+},j+2} \mbox{ \ for all } 1 \leqslant j \leqslant n-2 \mbox{ and } m^{+}_0 +j \mbox{ is odd.}
  \end{align*}
\end{proof}

\begin{definition}
  We refer to the bicentrally balanced $C_4$-face-magic projective labeling $Z$ in Theorem \ref{thmStandardProjectiveLabeling} as the
  \textit{standard projective labeling associated with} $X$.
  We say that $Z$ is a standard bicentrally balanced $C_4$-face-magic projective labeling on $\mathcal{P}_{m,n}$.
\end{definition}

As a result of Theorem \ref{thmStandardProjectiveLabeling}, we only need to find the
standard bicentrally balanced $C_4$-face-magic projective labelings on $\mathcal{P}_{m,n}$.

\begin{example}
Table~\ref{tableProj9x9Labelingexample} illustrates a standard bicentrally balanced $C_4$-face-magic projective
labeling on $\mathcal{P}_{9,9}$.
For convenience, we display the $9 \times 9$ projective grid graph as a $9 \times 9$ checkerboard where each square cell represents a vertex
  and square cells that share an edge are adjacent.
\end{example}

\begin{table}[h]
\hspace{0.0in}\begin{tabular}{|c|c|c|c|c|c|c|c|c|}
\hline
31  &51  &32   &47  &36  &46  &40  &42  &41    \\[0.05in]
\hline
53  &30  &52   &34  &48  &35  &44  &39  &43    \\[0.05in]
\hline
28  &54  &29   &50  &33  &49  &37  &45  &38    \\[0.05in]
\hline
65  &18  &64   &22  &60  &23  &56  &27  &55    \\[0.05in]
\hline
16   &66  &17   &62  &21  &61  &25  &57  &26     \\[0.05in]
\hline
68  &15  &67   &19  &63  &20  &59  &24  &58    \\[0.05in]
\hline
4  &78  &5   &74  &9   &73  &13  &69  &14     \\[0.05in]
\hline
80  &3   &79   &7   &75  &8   &71  &12  &70     \\[0.05in]
\hline
1   &81  &2    &77  &6   &76  &10  &72  &11     \\[0.05in]
\hline
\end{tabular}
\caption{A standard bicentrally balanced $C_4$-face magic projective labeling on $\mathcal{P}_{9,9}$ with $C_4$-face-magic value 165.}
\label{tableProj9x9Labelingexample}
\end{table}

\subsection{Partial bicentrally balanced labeling}

We need to introduce labelings on subgrids of $\mathcal{P}_{m,n}$ in order to determine a category of
standard bicentrally balanced $C_4$-face-magic labelings on $\mathcal{P}_{m,n}$.

\begin{definition} \label{defnPartialLabelings}
  Let $m\geqslant 3$ and $n \geqslant 3$ be odd integers, and let $M \geqslant m$ and $N \geqslant n$ be odd integers.
  Let $\mathrm{Grid}(m,n) =\{ (i,j) : 1 \leqslant i \leqslant m \mbox{ and } 1 \leqslant j \leqslant n \}$ be the
  $m \times n$ subgrid of $\mathcal{P}_{M,N}$.
  Let $P_m \times P_n$ represent the $m \times  n$ planar grid subgraph of $\mathcal{P}_{M,N}$ on $\mathrm{Grid}(m,n)$.
  Let $X= \{ x_{i,j} : (i,j)\in \mathrm{Grid}(m,n) \}$ be a labeling on $P_m \times P_n$
  We say $X$ is a \textit{partial bicentrally balanced $C_4$-face-magic labeling on the $m \times n$ subgrid of $\mathcal{P}_{M,N}$} if
  \begin{enumerate}
    \item the sum of the labels on the vertices of each $C_4$-face of $P_m \times P_n$ is $2MN+3$,
    \item $\{ x_{i,j} : (i,j) \in \mathrm{Grid}(m,n) \mbox{ and } i+j \mbox{ is even}\} = \{ 1,2,\ldots,\tfrac{1}{2}mn+\frac{1}{2}\}$,
    \item $\{ x_{i,j} : (i,j) \in \mathrm{Grid}(m,n) \mbox{ and } i+j \mbox{ is odd}\} = \{ MN- \tfrac{1}{2}mn + \tfrac{3}{2}, MN- \tfrac{1}{2}mn + \tfrac{5}{2},\ldots,MN\}$,
    \item $x_{i,j}+ x_{m+1-i,n+1-j} = \tfrac{1}{2} mn+ \tfrac{3}{2}$ if $(i,j)\in\mathrm{Grid}(m,n)$ and $i+j$ is even,
    \item $x_{i,j}+ x_{m+1-i,n+1-j} = 2MN - \tfrac{1}{2} mn + \tfrac{3}{2}$ if $(i,j)\in\mathrm{Grid}(m,n)$ and $i+j$ is odd, and
    \item $X$ satisfies
    \begin{align*}
    x_{i,n_0^{+}} &<  x_{i+2,n_0^{+}} \mbox{ \ for all } 1 \leqslant i \leqslant m-2 \mbox{ and } i+ n^{+}_0 \mbox{ is even,}  \\
    x_{i,n_0^{+}} &>  x_{i+2,n_0^{+}} \mbox{ \ for all } 1 \leqslant i \leqslant m-2 \mbox{ and } i+ n^{+}_0 \mbox{ is odd,}  \\
    x_{m_0^{+},j} &<  x_{m_0^{+},j+2} \mbox{ \ for all } 1 \leqslant j \leqslant n-2 \mbox{ and } m^{+}_0 +j \mbox{ is even, and} \\
    x_{m_0^{+},j} &>  x_{m_0^{+},j+2} \mbox{ \ for all } 1 \leqslant j \leqslant n-2 \mbox{ and } m^{+}_0 +j \mbox{ is odd.}
  \end{align*}
  \end{enumerate}
\end{definition}

\subsection{Partial alternating lexicographic labeling}

We introduce two partial bicentrally balanced $C_4$-face-magic labelings on an $m \times n$ subgrid of $\mathcal{P}_{M,N}$.

\begin{definition} \label{defnBaseAlternatingLexicographicLabeling}
  Let $m\geqslant 3$ and $n \geqslant 3$ be odd integers, and let $M \geqslant m$ and $N \geqslant n$ be odd integers.
  The \textit{partial horizontal alternating lexicographic labeling on the $m \times n$ subgrid of $\mathcal{P}_{M,N}$},
  denoted by $\mathrm{HALL}_{M,N}(m,n)$, is the labeling
  $\mathrm{HALL}_{M,N}(m,n)=\{x_{i,j} : (i,j)\in \mathrm{Grid}(m,n) \}$ given by
  \begin{itemize}
  \item $x_{2i-1,2j-1} = m(j-1)+i$, for all $1 \leqslant i \leqslant m^{+}_0$ and $1 \leqslant j \leqslant n^{+}_0$,
  \item $x_{2i,2j} = m(j-1)+ m^{+}_0 +i$, for all $1 \leqslant i \leqslant m_0$ and $1 \leqslant j \leqslant n_0$,
  \item $x_{2i,2j-1} = MN+ m(1-j)+ 1 -i$, for all $1 \leqslant i \leqslant m_0$ and $1 \leqslant j \leqslant n^{+}_0$, and
  \item $x_{2i-1,2j} = MN - mj+ m^{+}_0 +1 -i$, for all $1 \leqslant i \leqslant m^{+}_0$ and $1 \leqslant j \leqslant n_0$.
  \end{itemize}

  Similarly, the \textit{partial vertical alternating lexicographic labeling on the $m \times n$ subgrid of $\mathcal{P}_{M,N}$},
  denoted by $\mathrm{VALL}_{M,N}(m,n)$, is the labeling
  $\mathrm{VALL}_{M,N}(m,n)=\{y_{i,j} : (i,j)\in \mathrm{Grid}(m,n) \}$ given by
  \begin{itemize}
  \item $y_{2i-1,2j-1} = n(i-1)+j$, for all $1 \leqslant i \leqslant m^{+}_0$ and $1 \leqslant j \leqslant n^{+}_0$,
  \item $y_{2i,2j} = n(i-1)+ n^{+}_0 +j$, for all $1 \leqslant i \leqslant m_0$ and $1 \leqslant j \leqslant n_0$,
  \item $y_{2i,2j-1} = MN - ni + n^{+}_0 + 1 -j$, for all $1 \leqslant i \leqslant m_0$ and $1 \leqslant j \leqslant n^{+}_0$, and
  \item $y_{2i-1,2j} = MN + n(1-i) +1 -j$, for all $1 \leqslant i \leqslant m^{+}_0$ and $1 \leqslant j \leqslant n_0$.
  \end{itemize}
\end{definition}

\begin{figure}
\begin{picture}(225,180)(-25,-5)

\multiput(0,0)(40,0){5}{\circle*{5}}
\multiput(0,40)(40,0){5}{\circle*{5}}
\multiput(0,80)(40,0){5}{\circle*{5}}
\multiput(0,120)(40,0){5}{\circle*{5}}
\multiput(0,160)(40,0){5}{\circle*{5}}

\put(0,0){\line(1,0){160}}
\put(0,40){\line(1,0){160}}
\put(0,80){\line(1,0){160}}
\put(0,120){\line(1,0){160}}
\put(0,160){\line(1,0){160}}

\put(0,0){\line(0,1){160}}
\put(40,0){\line(0,1){160}}
\put(80,0){\line(0,1){160}}
\put(120,0){\line(0,1){160}}
\put(160,0){\line(0,1){160}}

\put(2,5){${\scriptstyle y_{1,1}=1}$}
\put(42,5){${\scriptstyle y_{2,1}=73}$}
\put(82,5){${\scriptstyle y_{3,1}=6}$}
\put(122,5){${\scriptstyle y_{4,1}=68}$}
\put(162,5){${\scriptstyle y_{5,1}=11}$}

\put(2,45){${\scriptstyle y_{1,2}=75}$}
\put(42,45){${\scriptstyle y_{2,2}=4}$}
\put(82,45){${\scriptstyle y_{3,2}=70}$}
\put(122,45){${\scriptstyle y_{4,2}=9}$}
\put(162,45){${\scriptstyle y_{5,2}=65}$}

\put(2,85){${\scriptstyle y_{1,3}=2}$}
\put(42,85){${\scriptstyle y_{2,3}=72}$}
\put(82,85){${\scriptstyle y_{3,3}=7}$}
\put(122,85){${\scriptstyle y_{4,3}=67}$}
\put(162,85){${\scriptstyle y_{5,3}=12}$}

\put(2,125){${\scriptstyle y_{1,4}=74}$}
\put(42,125){${\scriptstyle y_{2,4}=5}$}
\put(82,125){${\scriptstyle y_{3,4}=69}$}
\put(122,125){${\scriptstyle y_{4,4}=10}$}
\put(162,125){${\scriptstyle y_{5,4}=64}$}

\put(2,165){${\scriptstyle y_{1,5}=3}$}
\put(42,165){${\scriptstyle y_{2,5}=71}$}
\put(82,165){${\scriptstyle y_{3,5}=8}$}
\put(122,165){${\scriptstyle y_{4,5}=66}$}
\put(162,165){${\scriptstyle y_{5,5}=13}$}

\end{picture}
\caption{The partial bicentrally balanced $C_4$-face-magic projective labeling $\mathrm{VALL}_{15,5}(5,5)$.}
\label{figVALLFiveByFiveProjectiveGridGraphLabeling}
\end{figure}
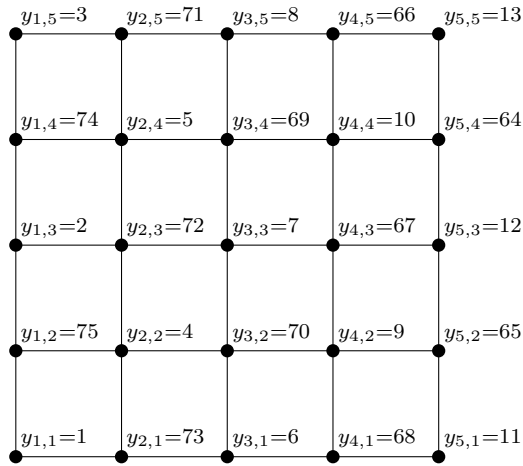

\begin{example}
  The labeling in Fig. \ref{figVALLFiveByFiveProjectiveGridGraphLabeling} is the partial vertical
  alternating lexicographic labeling on the $5 \times 5$ subgrid of $\mathcal{P}_{15,5}$ denoted by $\mathrm{VALL}_{15,5}(5,5)$.
\end{example}

\begin{theorem} \label{thmHALLandVALLProjectiveLabelings}
  Let $m\geqslant 3$ and $n \geqslant 3$ be odd integers, and let $M \geqslant m$ and $N \geqslant n$ be odd integers.
  The horizontal alternating lexicographic labeling $\mathrm{HALL}_{M,N}(m,n)$ and the vertical alternating lexicographic labeling $\mathrm{VALL}_{M,N}(m,n)$
  are partial bicentrally balanced $C_4$-face-magic labelings on the $m \times n$ subgrid of $\mathcal{P}_{M,N}$.
\end{theorem}

\begin{proof}
  We show that $\mathrm{VALL}_{M,N}(m,n)$ is a partial bicentrally balanced $C_4$-face-magic labelings on the $m \times n$ subgrid of $\mathcal{P}_{M,N}$.
  The proof that $\mathrm{HALL}_{M,N}(m,n)$ is a partial bicentrally balanced $C_4$-face-magic labelings on the $m \times n$ subgrid of $\mathcal{P}_{M,N}$ is similar.

  We observe that for the vertices $(i,j)$ where $i+j$ even, we assign the labels $1,2,\ldots,\tfrac{1}{2}mn+\tfrac{1}{2}$ in lexicographic order;
however, for the vertices $(i,j)$ where $i+j$ odd, we assign the labels $MN-(\tfrac{1}{2}mn-\tfrac{3}{2}),MN-(\tfrac{1}{2}mn-\tfrac{5}{2}),\ldots,MN$ in reverse lexicographic order.

We have  $y_{2i-1,2j-1} + y_{2i-1,2j}=MN+1$ for $1 \leqslant i \leqslant m_0+1$ and  $1 \leqslant j \leqslant n_0$.
Also, we have $y_{2i,2j-1} + y_{2i,2j}=MN+2$  for $1 \leqslant i \leqslant m_0$ and $1 \leqslant j \leqslant n_0$.
Thus, for $1 \leqslant i \leqslant m-1$ and  $1 \leqslant j \leqslant n_0$, we have
\begin{equation*}
  y_{i,2j-1} + y_{i,2j} + y_{i+1,2j-1} + y_{i+1,2j}=2MN+3.
\end{equation*}

Next, we have $y_{2i-1,2j} + y_{2i-1,2j+1}=MN+2$ for $1 \leqslant i \leqslant m_0+1$ and $1 \leqslant j \leqslant n_0$.
Also, we have $y_{2i,2j} + y_{2i,2j+1}=MN+1$ for $1 \leqslant i \leqslant m_0$ and $1 \leqslant j \leqslant n_0$.
Thus, for $1 \leqslant i \leqslant m-1$ and  $1 \leqslant j \leqslant n_0$, we have
\begin{equation*}
  y_{i,2j} + y_{i,2j+1} + y_{i+1,2j} + y_{i+1,2j+1}=2MN+3.
\end{equation*}

We observe that, for $1\leqslant j \leqslant n_0 +1$, $y_{1,2j-1} +y_{m,n+2-2j}=\tfrac{1}{2}mn+\tfrac{3}{2}$, and,
for $1\leqslant j \leqslant n_0$, $y_{1,2j} +y_{m,n+1-2j}=2MN - \tfrac{1}{2}mn+\tfrac{3}{2}$.
Thus, for $1\leqslant j \leqslant n-1$, we have
\begin{equation*}
  y_{1,j} +y_{m,n+1-j} +  y_{1,j+1} +y_{m,n-j} =2MN+3.
\end{equation*}

Similarly, for $1\leqslant i \leqslant m_0 +1$, $y_{2i-1,1} +y_{m+2-2i,n}=\tfrac{1}{2}mn+\tfrac{3}{2}$, and,
for $1\leqslant i \leqslant m_0$, $y_{2i,1} +y_{m+1-2i,n}=2MN- \tfrac{1}{2}mn+\tfrac{3}{2}$.
Thus, for $1\leqslant i \leqslant m-1$, we have
\begin{equation*}
  y_{i,1} +y_{m+1-i,n} +  y_{i+1,1} +y_{m-i,n} =2MN+3.
\end{equation*}

For all $1 \leqslant i \leqslant m^{+}_0$ and $1 \leqslant j \leqslant n^{+}_0$,
\begin{equation*}
  y_{2i-1,2j-1} + y_{m+2-2i,n+2-2j} = \tfrac{1}{2} mn + \tfrac{3}{2}.
\end{equation*}
Also, for all $1 \leqslant i \leqslant m_0$ and $1 \leqslant j \leqslant n_0$,
\begin{equation*}
  y_{2i,2j} + y_{m+1-2i,n+1-2j} = \tfrac{1}{2} mn + \tfrac{3}{2}.
\end{equation*}
We observe that, for all $1 \leqslant i \leqslant m^{+}_0$ and $1 \leqslant j \leqslant n_0$,
\begin{equation*}
  y_{2i-1,2j} + y_{m+2-2i,n+1-2j} = 2MN -\tfrac{1}{2} mn + \tfrac{3}{2}.
\end{equation*}
Similarly, for all $1 \leqslant i \leqslant m_0$ and $1 \leqslant j \leqslant n^{+}_0$,
\begin{equation*}
  y_{2i,2j-1} + y_{m+1-2i,n+2-2j} = 2MN -\tfrac{1}{2} mn + \tfrac{3}{2}.
\end{equation*}

Suppose $n\equiv 1\, (\mathrm{mod}\, 4)$. Then $n^{+}_0$ is odd, and we let $n'_0$ be the  positive integer such that $n^{+}_0=2n'_0 -1$.
Thus
\begin{align*}
  y_{2i-1,2n'_0 -1} &= n(i-1)+n'_0   \mbox{ \ for all }  1 \leqslant i \leqslant m^{+}_0, \mbox{ and } \\
  y_{2i,2n'_0 -1} &= MN - ni + n'_0    \mbox{ \ for all  } 1 \leqslant i \leqslant m_0.
\end{align*}
Hence,
\begin{align}
  y_{i,n_0^{+}} &<  y_{i+2,n_0^{+}} \mbox{ \ for all } 1 \leqslant i \leqslant m-2 \mbox{ and } i+ n^{+}_0 \mbox{ is even, and} \label{eqnStandardCentralColumnIncrease} \\
  y_{i,n_0^{+}} &>  y_{i+2,n_0^{+}} \mbox{ \ for all } 1 \leqslant i \leqslant m-2 \mbox{ and } i+ n^{+}_0 \mbox{ is odd.} \label{eqnStandardCentralColumnDecrease}
  \end{align}
A similar argument shows that (\ref{eqnStandardCentralColumnIncrease}) and (\ref{eqnStandardCentralColumnDecrease}) hold when $n\equiv 3\, (\mathrm{mod}\, 4)$.
Also, a similar argument shows that
    \begin{align*}
    y_{m_0^{+},j} &<  y_{m_0^{+},j+2} \mbox{ \ for all } 1 \leqslant j \leqslant n-2 \mbox{ and } m^{+}_0 +j \mbox{ is even, and} \\
    y_{m_0^{+},j} &>  y_{m_0^{+},j+2} \mbox{ \ for all } 1 \leqslant j \leqslant n-2 \mbox{ and } m^{+}_0 +j \mbox{ is odd.}
  \end{align*}
\end{proof}

\subsection{Alternating connected sums}

\begin{definition} \label{defnHorizontalAlternatingConnectedSum}
  Let $r \geqslant 3$ be an odd integer and let $r_0$ be the positive integer such that $r=2r_0 +1$.
  Let $X$ be a partial  bicentrally balanced $C_4$-face-magic labeling on the $m \times n$ subgrid of $\mathcal{P}_{M,N}$.
  Suppose $r \leqslant M/m$.
  The \textit{$r$-horizontal alternating connected sum of $X$}, denoted by $\mathrm{HASC}^r(X)$,
  is the labeling $Y=\{ y_{i,j} : (i,j) \in \mathrm{Grid}(mr,n) \}$ on $\mathrm{Grid}(mr,n)$ given by,
  for all $(i,j) \in \mathrm{Grid}(m,n)$ and for all integers $0 \leqslant k \leqslant r_0$,
  \begin{itemize}
    \item $ y_{(2k-1)m+i,j} = MN - x_{i,j} +(k-\tfrac{1}{2})(mn) +\tfrac{3}{2}$ if $1 \leqslant k \leqslant r_0$ and $i+j$ is odd,
    \item $ y_{(2k-1)m+i,j} = MN - x_{i,j} -(k-\tfrac{1}{2})(mn) + \tfrac{3}{2}$ if $1 \leqslant k \leqslant r_0$ and $i+j$ is even,
    \item $ y_{(2k)m+i,j} = x_{i,j} -k(mn)$ if $0 \leqslant k \leqslant r_0$ and $i+j$ is odd, and
    \item $ y_{(2k)m+i,j} = x_{i,j} +k(mn)$ if $0 \leqslant k \leqslant r_0$ and $i+j$ is even.
  \end{itemize}

  Suppose $r \leqslant N/n$.
  The \textit{$r$-vertical alternating connected sum of $X$}, denoted by $\mathrm{VASC}^r(X)$,
  is the labeling $Y=\{ y_{i,j} : (i,j) \in \mathrm{Grid}(m,nr) \}$ on $\mathrm{Grid}(m,nr)$ given by,
  for all $(i,j) \in \mathrm{Grid}(m,n)$ and for all integers $0 \leqslant k \leqslant r_0$,
  \begin{itemize}
    \item $ y_{i,(2k-1)n+j} = MN - x_{i,j} +(k-\tfrac{1}{2})(mn) +\tfrac{3}{2}$ if $1 \leqslant k \leqslant r_0$ and $i+j$ is odd,
    \item $ y_{i,(2k-1)n+j} = MN - x_{i,j} -(k-\tfrac{1}{2})(mn)  + \tfrac{3}{2}$ if $1 \leqslant k \leqslant r_0$ and $i+j$ is even,
    \item $ y_{i,(2k)n+j} = x_{i,j} -k(mn)$ if $0 \leqslant k \leqslant r_0$ and $i+j$ is odd, and
    \item $ y_{i,(2k)n+j} = x_{i,j} +k(mn)$ if $0 \leqslant k \leqslant r_0$ and $i+j$ is even.
  \end{itemize}
\end{definition}

\begin{example}
  The 3-horizontal alternating connected sum of $\mathrm{HALL}_{15,5}(5,5)$ is given in Table \ref{tableHorizontalAltConnSumofHALLexample}.
  For convenience, we display the $15 \times 5$ projective grid graph as a $15 \times 5$ checkerboard.
\end{example}

\begin{table}[h]
\hspace{0.0in}\begin{tabular}{|c|c|c|c|c|c|c|c|c|c|c|c|c|c|c|}
\hline
11  &65  &12   &64  &13  &53  &24  &52  &25   &51  &36  &40  &37  &39  &38   \\[0.05in]
\hline
68  &9   &67   &10  &66  &21  &55  &22  &54   &23  &43  &34  &42  &35  &41  \\[0.05in]
\hline
6   &70  &7    &69  &8   &58  &19  &57  &20   &56  &31  &45  &32  &44  &33   \\[0.05in]
\hline
73  &4   &72   &5   &71  &16  &60  &17  &59   &18  &48  &29  &47  &30  &46    \\[0.05in]
\hline
1    &75  &2   &74  &3   &63  &14  &62  &15   &61  &26  &50  &27  &49  &28    \\[0.05in]
\hline
\end{tabular}
\caption{A 3-horizontal alternating connected sum of $\mathrm{HALL}_{15,5}(5,5)$ on $\mathcal{P}_{15,5}$ with $C_4$-face-magic value 153.}
\label{tableHorizontalAltConnSumofHALLexample}
\end{table}

\begin{theorem} \label{thmHASCandVASCareBicentrallyBalanced}
  Suppose $X$ is a partial  bicentrally balanced $C_4$-face-magic labeling on the $m \times n$ subgrid of $\mathcal{P}_{M,N}$.
  If $r$ is an odd positive integer such that $mr\leqslant M$, then
  the $r$-horizontal alternating connected sum of $X$, $\mathrm{HASC}^r(X)$, is a partial  bicentrally balanced $C_4$-face-magic labelings on the $mr \times n$ subgrid of $\mathcal{P}_{M,N}$.
  Similarly, if $r$ is an odd positive integer such that $nr\leqslant N$, then
  the $r$-vertical alternating connected sum of $X$, $\mathrm{VASC}^r(X)$, is a partial  bicentrally balanced $C_4$-face-magic labelings on the $m \times nr$ subgrid of $\mathcal{P}_{M,N}$.
\end{theorem}

\begin{proof}
  We show that $\mathrm{HASC}^r(X)$, is a partial  bicentrally balanced $C_4$-face-magic labeling on the $mr \times n$ subgrid of $\mathcal{P}_{M,N}$.
  The proof that $\mathrm{VASC}^r(X)$, is a partial  bicentrally balanced $C_4$-face-magic labeling on the $m \times nr$ subgrid of $\mathcal{P}_{M,N}$ is similar.

  First, when $1\leqslant k \leqslant r_0$, we have
  \begin{align*}
    y_{(2k-1)m+2i-1,j} + y_{(2k-1)m+2i-1,j+1}   &=   2MN+3 -( x_{2i-1,j} + x_{2i-1,j+1} ) \mbox{ and } \\
    y_{(2k-1)m+2i,j} + y_{(2k-1)m+2i,j+1}       &=    2MN+3 -( x_{2i,j} + x_{2i,j+1} ).
  \end{align*}
  Thus, for all $1\leqslant i \leqslant m-1$, $1\leqslant j \leqslant n-1$, and $1\leqslant k \leqslant r_0$, we have
  \begin{equation*}
    y_{(2k-1)m+i,j} + y_{(2k-1)m+i,j+1} + y_{(2k-1)m+i+1,j} + y_{(2k-1)m+i+1,j+1} = 2MN+3.
  \end{equation*}

  Second, when $0\leqslant k \leqslant r_0$, we have
  \begin{align*}
    y_{(2k)m+2i-1,j} + y_{(2k)m+2i-1,j+1}   &=   x_{2i-1,j} + x_{2i-1,j+1}  \mbox{ and } \\
    y_{(2k)m+2i,j} + y_{(2k)m+2i,j+1}       &=    x_{2i,j} + x_{2i,j+1}.
  \end{align*}
  Thus, for all $1\leqslant i \leqslant m-1$, $1\leqslant j \leqslant n-1$, and $0\leqslant k \leqslant r_0$, we have
  \begin{equation*}
    y_{(2k)m+i,j} + y_{(2k)m+i,j+1} + y_{(2k)m+i+1,j} + y_{(2k)m+i+1,j+1} = 2MN+3.
  \end{equation*}

  When we set the $C_4$-face sums below equal,
  \begin{equation*}
    x_{i,j} + x_{i,j+1} + x_{i+1,j} + x_{i+1,j+1} = x_{i+1,j} + x_{i+1,j+1} + x_{i+2,j} + x_{i+2,j+1},
  \end{equation*}
  we obtain
  \begin{equation*}
    x_{i,j} + x_{i,j+1}  =  x_{i+2,j} + x_{i+2,j+1}.
  \end{equation*}
  Thus, for all $1 \leqslant j \leqslant n-1$, we have
  \begin{equation} \label{eqnOneSideSumEqualsOtherSide}
    x_{1,j} + x_{1,j+1}  =  x_{m,j} + x_{m,j+1}.
  \end{equation}

  Third, for all $1 \leqslant j \leqslant n-1$ and $1 \leqslant k \leqslant r_0$, we have
  \begin{equation*}
     y_{(2k)m,j} + y_{(2k)m,j+1}  = 2MN+3 -( x_{m,j} + x_{m,j+1} ).
  \end{equation*}
  Similarly, we have
  \begin{equation*}
     y_{(2k)m+1,j} + y_{(2k)m+1,j+1}  =  x_{1,j} + x_{1,j+1}.
  \end{equation*}
  By (\ref{eqnOneSideSumEqualsOtherSide}), we have
  \begin{equation*}
      y_{(2k)m,j} + y_{(2k)m,j+1} + y_{(2k)m+1,j} + y_{(2k)m+1,j+1}  =  2MN+3.
  \end{equation*}
  A similar argument shows that, for all $1 \leqslant j \leqslant n-1$ and $1 \leqslant k \leqslant r_0$, we have
  \begin{equation*}
      y_{(2k-1)m,j} + y_{(2k-1)m,j+1} + y_{(2k-1)m+1,j} + y_{(2k-1)m+1,j+1}  =  2MN+3.
  \end{equation*}

  From the construction of $Y$, we have
  \begin{align*}
    \{ y_{i,j} : (i,j) \in \mathrm{Grid}(mr,n) \mbox{ and } i+j \mbox{ is even}\} &= \{ 1,2,\ldots,\tfrac{1}{2}mnr+\tfrac{1}{2}\}, \mbox{ and } \\
    \{ y_{i,j} : (i,j) \in \mathrm{Grid}(mr,n) \mbox{ and } i+j \mbox{ is odd}\} &= \{ MN- \tfrac{1}{2}mnr + \tfrac{3}{2}, \\
    &MN- \tfrac{1}{2}mnr + \tfrac{5}{2}, \allowbreak \ldots,MN\}.
  \end{align*}
   If $i+j$ is even, we have
   \begin{align*}
     y_{(2k-1)m+i,j} + y_{(2(r_0-k+1)-1)m +(m+1-i),n+1-j} &= 2MN -\tfrac{1}{2} mnr +\tfrac{3}{2}, \mbox{ \ and } \\
     y_{(2k)m+i,j} + y_{2(r_0-k)m +(m+1-i),n+1-j}  &= \tfrac{1}{2} mnr +\tfrac{3}{2}.
   \end{align*}
   If $i+j$ is odd, we have
   \begin{align*}
     y_{(2k-1)m+i,j} + y_{(2(r_0-k+1)-1)m +(m+1-i),n+1-j}  &= \tfrac{1}{2} mnr +\tfrac{3}{2}, \mbox{ \ and }  \\
      y_{(2k)m+i,j} + y_{2(r_0-k)m +(m+1-i),n+1-j} &= 2MN -\tfrac{1}{2} mnr +\tfrac{3}{2}.
   \end{align*}
   Thus, for all $1 \leqslant i \leqslant mr$ and $1 \leqslant j \leqslant n$,
   \begin{align*}
     y_{i,j} + y_{mr+1-i,n+1-j} &= \tfrac{1}{2} mnr +\tfrac{3}{2}, \mbox{ \ if } i+j \mbox{ is even, and} \\
     y_{i,j} + y_{mr+1-i,n+1-j} &= 2MN -\tfrac{1}{2} mnr +\tfrac{3}{2}, \mbox{ \ if } i+j \mbox{ is odd.}
   \end{align*}

   We have,
   for all $1 \leqslant i \leqslant m$ and $0 \leqslant k \leqslant r_0$,
  \begin{align*}
   y_{(2k-1)m+i,n^{+}_0} &= MN - x_{i,n^{+}_0} +(k-\tfrac{1}{2})(mn) +\tfrac{3}{2} \mbox{ \ if } i+n^{+}_0 \mbox{ is odd and} \\
   y_{(2k)m+i,n^{+}_0} &= x_{i,n^{+}_0} +k(mn) \mbox{ \ if } i+n^{+}_0 \mbox{ is even.}
  \end{align*}
   Since, for all  $1 \leqslant i \leqslant m-2$,
   \begin{align*}
    x_{i,n_0^{+}} &<  x_{i+2,n_0^{+}} \mbox{ \ if } i+ n^{+}_0 \mbox{ is even and }    \\
    x_{i,n_0^{+}} &>  x_{i+2,n_0^{+}} \mbox{ \ if } i+ n^{+}_0 \mbox{ is odd,}
  \end{align*}
   we have
   \begin{align*}
    y_{(2k-1)m+i,n_0^{+}} &<  y_{(2k-1)m+i+2,n_0^{+}} \mbox{ \ if }  (2k-1)m+i+ n^{+}_0 \mbox{ is even and }  \\
    y_{(2k)m+i,n_0^{+}} &<  y_{(2k)m+i+2,n_0^{+}} \mbox{ \ if } (2k)m+i+ n^{+}_0 \mbox{ is even.}
  \end{align*}
   In addition, we have
   \begin{align*}
     (k-1)(mn) +\tfrac{1}{2}mn +\tfrac{3}{2} &\leqslant y_{(2k-1)m+i,n^{+}_0} \leqslant k(mn) \\
     &\mbox{ \ if } (2k-1)m+i+n^{+}_0 \mbox{ is even and } \\
     k(mn) + 1 &\leqslant y_{(2k)m+i,n^{+}_0} \leqslant k(mn) +\tfrac{1}{2}mn +\tfrac{1}{2}  \\
     &\mbox{ \ if } (2k)m+i+n^{+}_0 \mbox{ is even.}
   \end{align*}
   Hence, for all $1 \leqslant i \leqslant mr-2$,
   \begin{equation*}
      y_{i,n_0^{+}} <  y_{i+2,n_0^{+}} \mbox{ \ if } i+ n^{+}_0 \mbox{ is even.}
   \end{equation*}
   A similar argument shows that, for all $1 \leqslant i \leqslant mr-2$,
   \begin{equation*}
      y_{i,n_0^{+}} >  y_{i+2,n_0^{+}} \mbox{ \ if } i+ n^{+}_0 \mbox{ is odd.}
   \end{equation*}

   Let $M^{+}_0$ be the positive integer such that $mr= 2M^{+}_0 -1$.
   Then $M^{+}_0 = m r_0 + m^{+}_0$.
   Suppose $r_0$ is odd. Let $r'_0$ be the positive integer such that $r_0 = 2 r'_0 -1$.
   For all $1 \leqslant j \leqslant n$, we have
   \begin{align*}
      y_{(2r'_0 -1)m+m^{+}_0,j} &= MN - x_{m^{+}_0,j} +(r'_0 -\tfrac{1}{2})(mn)  +\tfrac{3}{2} \mbox{ \ if }  m^{+}_0+j \mbox{ is odd and } \\
      y_{(2r'_0 -1)m+m^{+}_0,j} &= MN - x_{m^{+}_0,j} - (r'_0 -\tfrac{1}{2})(mn)  +\tfrac{3}{2} \mbox{ \ if }  m^{+}_0+j \mbox{ is even.}
   \end{align*}
   Since, for all  $1 \leqslant j \leqslant n-2$,
   \begin{align*}
    x_{m_0^{+},j} &<  x_{m_0^{+},j+2} \mbox{ \ if } m^{+}_0 +j  \mbox{ is even and }  \\
    x_{m_0^{+},j} &>  x_{m_0^{+},j+2} \mbox{ \ if } m^{+}_0 +j  \mbox{ is odd,}
  \end{align*}
   we have
   \begin{align*}
    y_{(2r'_0 -1)m+ m_0^{+},j} &<  y_{(2r'_0 -1)m+ m_0^{+},j+2} \mbox{ \ if }  (2r'_0 -1)m+ m^{+}_0 +j \mbox{ is even and }  \\
    y_{(2r'_0 -1)m+ m_0^{+},j} &>  y_{(2r'_0 -1)m+ m_0^{+},j+2} \mbox{ \ if }  (2r'_0 -1)m+ m^{+}_0 +j \mbox{ is odd.}
  \end{align*}
   Thus, for all $1 \leqslant j \leqslant n-2$, we have
   \begin{align}
    y_{M_0^{+},j} &<  y_{M_0^{+},j+2}   \mbox{ \ if } M^{+}_0 +j \mbox{ is even and } \label{eqnConnectedSumLabelCenterColumnIncrease2}  \\
    y_{M_0^{+},j} &>  y_{M_0^{+},j+2}   \mbox{ \ if } M^{+}_0 +j \mbox{ is odd.} \label{eqnConnectedSumLabelCenterColumnDecrease2}
  \end{align}
   A similar argument shows that (\ref{eqnConnectedSumLabelCenterColumnIncrease2}) and (\ref{eqnConnectedSumLabelCenterColumnDecrease2}) hold when $r_0$ is even.
\end{proof}

\subsection{Labelings associated with a projective factorization sequence}

\begin{definition} \label{defnBicentrallyBalancedLabelingAssociatedWithF}
  Let $m\geqslant 3$ and $n\geqslant 3$ be odd integers.
   \begin{enumerate}
     \item Let $F= (m_i, n_i : 1 \leqslant i \leqslant k)$ be an $(m,n)$-projective factorization sequence.
     See Definition \ref{defnProjectiveFactorizationSequence}.
       Let $X_1 =\mathrm{HALL}(m_1,n_1)$.
       For $2 \leqslant i \leqslant k$, let $Y_i =\mathrm{HACS}^{m_i}(X_{i-1})$ and $X_i =\mathrm{VACS}^{n_i}(Y_{i})$.
       The \textit{horizontal bicentrally balanced labeling associated with} $F$ is denoted by $\mathrm{HBBL}(F)=X_k$.
      \item Let $F'= (n'_i, m'_i : 1 \leqslant i \leqslant k)$ be an $(n,m)$-projective factorization sequence.
       Let $X'_1 =\mathrm{VALL}(m'_1,n'_1)$.
       For $2 \leqslant i \leqslant k$, let $Y'_i =\mathrm{VACS}^{n'_i}(X'_{i-1})$ and $X'_i =\mathrm{HACS}^{m'_i}(Y'_{i})$.
       The \textit{vertical bicentrally balanced labeling associated with} $F'$ is denoted by $\mathrm{VBBL}(F')=X'_k$.
   \end{enumerate}
\end{definition}

\begin{theorem} \label{thmStandardProjectiveC4MagicBicentralConstruction}
  Let $m\geqslant 3$ and $n\geqslant 3$ be odd integers.
   Suppose $X$ is constructed in one of the following two ways.
   \begin{enumerate}
     \item Suppose $X= \mathrm{HBBL}(F)$ for some $(m,n)$-projective factorization sequence
     $F=(m_1,\allowbreak n_1,\allowbreak m_2,\allowbreak n_2,\ldots,m_{k},n_{k})$.
     \item Suppose $X= \mathrm{VBBL}(F')$ for some $(n,m)$-projective factorization sequence
     $F'=(n'_1,\allowbreak m'_1,\allowbreak n'_2,\allowbreak m'_2,\ldots,n'_{k},m'_{k})$.
   \end{enumerate}
    Then $X$ is a standard bicentrally balanced $C_4$-face-magic projective labeling on $\mathcal{P}_{m,n}$.

   Furthermore, distinct $(m,n)$-projective factorization sequences $F_1$ and $F_2$ give rise
   to distinct standard bicentrally balanced $C_4$-face-magic projective labelings $\mathrm{HBBL} \allowbreak (F_1)$ and $\mathrm{HBBL} \allowbreak (F_2)$ on $\mathcal{P}_{m,n}$.
   Similarly, distinct $(n,m)$-projective factorization sequences $F'_1$ and $F'_2$ give rise
   to distinct standard bicentrally balanced $C_4$-face-magic projective labelings $\mathrm{VBBL}(F'_1)$ and $\mathrm{VBBL}(F'_2)$ on $\mathcal{P}_{m,n}$.
\end{theorem}

\begin{proof}
  We first show that $\mathrm{HBBL}(F)$ is a standard bicentrally balanced $C_4$-face-magic projective labeling on $\mathcal{P}_{m,n}$.
  By Theorem \ref{thmHALLandVALLProjectiveLabelings}, $X_1 =\mathrm{HALL}_{m,n}(m_1,n_1)$ is a partial
  bicentrally balanced $C_4$-face-magic labeling on the $m_1 \times n_1$ subgrid of $\mathcal{P}_{m,n}$.
  For $1 \leqslant i \leqslant k$, let
  \begin{align*}
    M_i &= m_1 m_2 \cdots m_i \mbox{ and } \\
    N_i &= n_1 n_2 \cdots n_i.
  \end{align*}
  For some integer $2 \leqslant i \leqslant k$, suppose $X_{i-1}$ is a  partial
  bicentrally balanced $C_4$-face-magic labeling on the $M_{i-1} \times N_{i-1}$ subgrid of $\mathcal{P}_{m,n}$.
  By Theorem \ref{thmHASCandVASCareBicentrallyBalanced}, $Y_i =\mathrm{HACS}^{m_i}(X_{i-1})$
  is a  partial bicentrally balanced $C_4$-face-magic labeling on the $M_{i} \times N_{i-1}$ subgrid of $\mathcal{P}_{m,n}$
  and $X_i =\mathrm{VACS}^{n_i}(Y_{i})$ is a  partial
  bicentrally balanced $C_4$-face-magic labeling on the $M_{i} \times N_{i}$ subgrid of $\mathcal{P}_{m,n}$
  Thus $X_k= \mathrm{HBBL}(F)$ is a  partial
  bicentrally balanced $C_4$-face-magic labeling on the $m \times n$ subgrid of $\mathcal{P}_{m,n}$.
  Hence,  $\mathrm{HBBL}(F)$ is a  standard
  bicentrally balanced $C_4$-face-magic projective labeling on $\mathcal{P}_{m,n}$.

  A similar argument shows that, for an $(n,m)$-projective factorization sequence $F'$,
  $\mathrm{VBBL}\allowbreak (F')$ is a standard bicentrally balanced $C_4$-face-magic projective labeling on $\mathcal{P}_{m,n}$.

  Let $F_i= (m_{i,j}, n_{i,j} : 1 \leqslant j \leqslant k_i)$, for $i=1,2$, be distinct $(m,n)$-projective factorization sequences.
  We need to show that $\mathrm{HBBL}(F_1)$ and $\mathrm{HBBL}(F_2)$ are distinct
  standard bicentrally balanced $C_4$-face-magic projective labelings on $\mathcal{P}_{m,n}$.
  Suppose $m_{1,j'}$ and $m_{2,j'}$ are the first entries in $F_1$ and $F_2$, respectively, such that  $m_{1,j'} \ne m_{2,j'}$.

  First, assume $j'=1$. Without loss of generality, we may assume $m_{1,1} < m_{2,1}$.
  Then the labels on $\mathrm{HBBL}(F_1)$ and $\mathrm{HBBL}(F_2)$ are the same on the $m_{1,1} \times 1$ subgrid of $\mathcal{P}_{m,n}$.
  Let $Z=\{y_{i,1} : 1 \leqslant i \leqslant m_{1,1}\}$ be the common labels of
  $\mathrm{HBBL}(F_1)$ and $\mathrm{HBBL}(F_2)$ on the $m_{1,1} \times 1$ subgrid of $\mathcal{P}_{m,n}$.
  Let $z$ be the smallest positive integer such that $z\in \mathrm{Label}_Z(\mathrm{Grid}(m_{1,1}, 1))$
  and $z+1\notin \mathrm{Label}_Z(\mathrm{Grid}(m_{1,1}, 1))$.
  Let $m'=m_{1,1}$.
  Then $z=y_{m',1}=\tfrac{1}{2}(m'+1)$.
  In $\mathrm{HBBL}(F_1)$ we have $y_{2,2}=z+1$, and
  in $\mathrm{HBBL}(F_2)$ we have $y_{m' +2,1}=z+1$.
  Thus $\mathrm{HBBL}(F_1)$ and $\mathrm{HBBL}(F_2)$ are distinct standard bicentrally balanced $C_4$-face-magic projective labelings on $\mathcal{P}_{m,n}$.

  Suppose $j'>1$.
  Without loss of generality, we may assume $m_{1,j'} < m_{2,j'}$.
  Let $m'=m_{1,j'} M_{j' -1}$ and $n'=N_{j' -1}$.
  Let $Z= \{ y_{i,j} : (i,j) \in \mathrm{Grid}(m',n')\}$ be the common labels of
  $\mathrm{HBBL}(F_1)$ and $\mathrm{HBBL}(F_2)$ on the $m' \times n'$ subgrid of $\mathcal{P}_{m,n}$.
  Let $z$ be the smallest positive integer such that $z\in \mathrm{Label}_{Z}(\mathrm{Grid}(m',n'))$
  and $z+1\notin \mathrm{Label}_{Z}(\mathrm{Grid}(m',n'))$.
  Then $z=y_{m',n'}=\tfrac{1}{2}(m' n'+1)$.
  In $\mathrm{HBBL}(F_1)$ we have $y_{2,n' +1}=z+1$ , and
  in $\mathrm{HBBL}(F_2)$ we have $y_{m' +2,1}=z+1$.
  Thus $\mathrm{HBBL}(F_1)$ and $\mathrm{HBBL}(F_2)$ are distinct standard bicentrally balanced $C_4$-face-magic projective labelings on $\mathcal{P}_{m,n}$.

  Suppose $n_{1,j'}$ and $n_{2,j'}$ are the first entries in $F_1$ and $F_2$, respectively, such that  $n_{1,j'} \ne n_{2,j'}$.
  A similar argument to the one above shows that $\mathrm{HBBL}(F_1)$ and $\mathrm{HBBL}(F_2)$ are distinct
  standard bicentrally balanced $C_4$-face-magic projective labelings on $\mathcal{P}_{m,n}$.

  Similarly, we can show that distinct $(n,m)$-projective factorization sequences $F'_1$ and $F'_2$ give rise
   to distinct standard bicentrally balanced $C_4$-face-magic projective labelings $\mathrm{VBBL}(F'_1)$ and $\mathrm{VBBL}(F'_2)$ on $\mathcal{P}_{m,n}$.
\end{proof}

\section{Enumerating bicentrally balanced labelings}

We will enumerate the minimum number of distinct $C_4$-face-magic projective labelings on $\mathcal{P}_{m,n}$.

\begin{notation}
  Let $m\geqslant 3$ be an odd integer.
  We define the function $\beta$ given by
  \begin{equation*}\label{eqnBetaFunction}
    \beta(m) = \left\{ \begin{array}{ll}
   \bigl( ( \tfrac{m-1}{4} )! \bigr)^2, & \mbox{ if } m\equiv 1\pmod{4}, \\[\medskipamount]
 ( \tfrac{m-3}{4} )! ( \tfrac{m+1}{4} )!, & \mbox{ if } m\equiv 3\pmod{4}.
\end{array} \right.
  \end{equation*}
\end{notation}

The following theorem gives us a lower bound on the number of distinct
$C_4$-face-magic projective labelings on $\mathcal{P}_{m,n}$ having $C_4$-face-magic value $2mn +1$ (or $2mn +3$) for distinct odd integers $m$ and $n$.

\begin{theorem}   \label{thmNumberOfBicentrallyBalancedLabelingsOnRectangle}
  Let $m \geqslant 3$ and $n \geqslant 3$ be distinct odd integers.
  Then the number of distinct $C_4$-face-magic projective labelings on $\mathcal{P}_{m,n}$ having $C_4$-face-magic value $2mn+1$ (or $2mn+3$)
  (up to symmetries on the projective plane) is at least
  \begin{equation*}
     \bigl( \tau(m,n) + \tau(n,m) \bigr) 2^{m/2+n/2-3} \beta(m)  \beta(n),
  \end{equation*}
  where $\tau(m,n)$ is the number of distinct $(m,n)$-projective factorization sequences. See Definition \ref{defnProjectiveFactorizationSequence}.
\end{theorem}

\begin{proof}
  Let $X$ be a standard bicentrally balanced $C_4$-face-magic projective labeling on $\mathcal{P}_{m,n}$.
  Let $\eta$ be a permutation on $\{1,2,\ldots,m_0\}$ such that $\eta(i)\equiv i \pmod{2}$ for all $1 \leqslant i \leqslant m_0$
  and $\mathcal{E}_{\eta}(X)$ be the labeling given in Definition \ref{defnElementaryProjectivePermutationColumns}.
  By Lemma \ref{lemmaElementaryProjectivePermutationColumns}, there are $\beta(m)$ distinct bicentrally balanced $C_4$-face-magic projective labelings
  of type $\mathcal{E}_{\eta}(X)$ on $\mathcal{P}_{m,n}$.
  Let $\alpha:\{1,2,\ldots,m_0\}\rightarrow \{0,1\}$ and $\mathcal{E}_{\alpha}(X)$ be the labeling given in Definition \ref{defnElementaryProjectiveSwapColumns}.
  By Lemma \ref{lemmaElementaryProjectiveSwapColumns}, there are $2^{m_0}$ distinct bicentrally balanced $C_4$-face-magic projective labelings
  of type $\mathcal{E}_{\alpha}(X)$ on $\mathcal{P}_{m,n}$.
  Similarly, by Lemma \ref{lemmaElementaryProjectivePermutationRows}, there are $\beta(n)$
  distinct bicentrally balanced $C_4$-face-magic projective labelings on $\mathcal{P}_{m,n}$ associated with the elementary projective labeling operation given
  in Definition \ref{defnElementaryProjectivePermutationRows}.
  Also, by Lemma \ref{lemmaElementaryProjectiveSwapRows}, there are $2^{n_0}$
  distinct bicentrally balanced $C_4$-face-magic projective labelings on $\mathcal{P}_{m,n}$ associated with the elementary projective labeling operation given
  in Definition \ref{defnElementaryProjectiveSwapRows}.
  Thus, there are $\beta(m) 2^{m_0} \beta(n) 2^{n_0}$ distinct bicentrally balanced $C_4$-face-magic projective labelings
  on $\mathcal{P}_{m,n}$ that are projective labeling equivalent to $X$.
  For each of these labelings, there are four labelings that are projective labeling equivalent by one of the symmetries
  $R_0$, $R_{180}$, $H$, or $V$ of the projective plane.
  Hence, there are $\tfrac{1}{4} \beta(m) 2^{m_0} \beta(n) 2^{n_0} = 2^{m/2+n/2-3} \beta(m) \beta(n)$ distinct bicentrally balanced $C_4$-face-magic projective labelings
  on $\mathcal{P}_{m,n}$ that are projective labeling equivalent to $X$ up to symmetries of the projective plane.
  By Remark \ref{remBicentrallyBalancedHasValue2MNPlus3} and Lemma \ref{lemmaProjectiveGridSemiC4FaceValueBalance}, a $C_4$-face-magic projective labeling on $\mathcal{P}_{m,n}$ is bicentrally balanced if and only if its $C_4$-face-magic value is $2mn+3$.
  Therefore, there are $2^{m/2+n/2-3} \beta(m) \beta(n)$ distinct $C_4$-face-magic projective labelings
  on $\mathcal{P}_{m,n}$ with $C_4$-face-magic value $2mn+3$ that are projective labeling equivalent to $X$ up to symmetries of the projective plane.

  By Theorem \ref{thmStandardProjectiveC4MagicBicentralConstruction}, each $(m,n)$-projective factorization sequence
  $F$ and each $(n,m)$-projective factorization sequence $F'$ are associated with unique standard
  bicentrally balanced $C_4$-face-magic projective labelings on $\mathcal{P}_{m,n}$
  given by $\mathrm{HBBL}(F)$ and $\mathrm{VBBL}(F')$, respectively.
  Therefore, there are at least
  \begin{equation*}
    \bigl( \tau(m,n) + \tau(n,m) \bigr) 2^{m/2 +n/2-3}  \beta(m) \allowbreak \beta(n)
  \end{equation*}
  distinct $C_4$-face-magic projective labelings on $\mathcal{P}_{m,n}$ (up to symmetries on the projective plane) that have $C_4$-face-magic value $2mn+3$.

  By Remark \ref{remOrderComplementLabelingCorrespondence}, there are at least
  \begin{equation*}
    \bigl( \tau(m,n) + \tau(n,m) \bigr) 2^{m/2 +n/2-3}  \beta(m) \allowbreak \beta(n)
  \end{equation*}
  distinct $C_4$-face-magic projective labelings on $\mathcal{P}_{m,n}$ (up to symmetries on the projective plane) that have $C_4$-face-magic value $2mn+1$.
\end{proof}

The next proposition gives us a lower bound on the number of distinct
$C_4$-face-magic projective labelings on $\mathcal{P}_{m,m}$ having $C_4$-face-magic value $2m^2 +1$ (or $2m^2 +3$).

\begin{theorem}   \label{thmNumberOfBicentrallyBalancedLabelingsOnSquare}
  Let $m \geqslant 3$ be an odd integer.
  Then the number of distinct $C_4$-face-magic projective labelings on $\mathcal{P}_{m,m}$ having $C_4$-face-magic value $2m^2 +1$ (or $2m^2 +3$)
  (up to symmetries on the projective plane) is at least
  \begin{equation*}
     \tau(m,m) 2^{m-3} \bigl( \beta(m) \bigr)^2.
  \end{equation*}
\end{theorem}

The proof of Theorem \ref{thmNumberOfBicentrallyBalancedLabelingsOnSquare} is similar to that of
Theorem \ref{thmNumberOfBicentrallyBalancedLabelingsOnRectangle}.

\begin{theorem} \label{thmNumberOfC4FaceMagicLabelingsOnRectangle}
  Let $m \geqslant 3$ and $n \geqslant 3$ be distinct odd integers.
  Then the number of distinct $C_4$-face-magic projective labelings on $\mathcal{P}_{m,n}$
  (up to symmetries on the projective plane) is at least
  \begin{equation*}
    \bigl(  \tau(m,n) + \tau(n,m) \bigr)   2^{m/2 +n/2-3} \bigl( (\tfrac{m-1}{2})! (\tfrac{n-1}{2})!   + 2 \beta(m)  \beta(n) \bigr).
  \end{equation*}
\end{theorem}

\begin{proof}
   By Lemma \ref{lemmaOddByOddProjectiveGridC4FaceValue}, the $C_4$-face-magic value of a $C_4$-face-magic projective labeling on
   $\mathcal{P}_{m,n}$ is either $2mn+1$, $2mn+2$, or $2mn+3$.
   By Theorem \ref{thmNumberOfCentrallyBalancedC4FaceMagicLabelingsOnRectangle}, there are exactly
   \begin{equation*}
    ( \tau(m,n) + \tau(n,m) )2^{m/2 +n/2-3} (\tfrac{m-1}{2})! (\tfrac{n-1}{2})!
   \end{equation*}
   distinct $C_4$-face-magic projective labelings on $\mathcal{P}_{m,n}$
  (up to symmetries on the projective plane) with $C_4$-face-magic value $2mn+2$.
  By Theorem \ref{thmNumberOfBicentrallyBalancedLabelingsOnRectangle}, there are at least
   \begin{equation*}
    ( \tau(m,n) + \tau(n,m) ) 2^{m/2+n/2-3} \beta(m)  \beta(n)
   \end{equation*}
   distinct $C_4$-face-magic projective labelings on $\mathcal{P}_{m,n}$
  (up to symmetries on the projective plane) with $C_4$-face-magic value $2mn+1$ (or $2mn+3$).
  The result follows.
\end{proof}

\begin{theorem} \label{thmNumberOfC4FaceMagicLabelingsOnSquare}
  Let $m \geqslant 3$ be an odd integer.
  Then the number of distinct $C_4$-face-magic projective labelings on $\mathcal{P}_{m,m}$
  (up to symmetries on the projective plane) is at least
  \begin{equation*}
    \tau(m,m)  2^{m-3} \bigl(  \bigl( (\tfrac{m-1}{2})! \bigr)^2 + 2 \bigl( \beta(m) \bigr)^2 \bigr).
  \end{equation*}
\end{theorem}

  The proof Theorem \ref{thmNumberOfC4FaceMagicLabelingsOnSquare} is similar to that of Theorem \ref{thmNumberOfC4FaceMagicLabelingsOnRectangle}.

\section{Open Problems}

We conjecture that the labelings given in Theorem \ref{thmStandardProjectiveC4MagicBicentralConstruction} are the only standard bicentrally balanced
$C_4$-face-magic projective labelings on $\mathcal{P}_{m,n}$.

\begin{conjecture}
Let $m \geqslant 3$ and $n \geqslant 3$ be odd integers.
Suppose $X$ is a standard bicentrally balanced $C_4$-face-magic projective labeling on $\mathcal{P}_{m,n}$.
Then either
\begin{enumerate}
  \item there exists an $(m,n)$-projective factorization sequence $F$ such that $X=\mbox{HBBL}(F)$ or
  \item there exists an $(n,m)$-projective factorization sequence $F'$ such that $X=\mbox{VBBL}(F')$.
\end{enumerate}
\end{conjecture}

Suppose $m\geqslant 3$ and $n \geqslant 3$ are odd integers.
The characterization of $C_4$-face-magic projective labelings  on $\mathcal{P}_{m,n}$ with $C_4$-face-magic value $2mn+2$
is given in \cite{Curran2}.
It is natural to ask if there is a characterization of $C_4$-face-magic projective labelings on $\mathcal{P}_{m,n}$ when $m\geqslant 2$ and
$n\geqslant 2$ are even integers.
By Lemma \ref{lemmaEvenByEvenProjectiveGridC4FaceValue}, the $C_4$-face-magic value of such a labeling is $2mn+2$.

\begin{problem}
  Let $m\geqslant 2$ and $n\geqslant 2$ be even integers.
  Find a characterization of the $C_4$-face-magic projective labelings on $\mathcal{P}_{m,n}$.
\end{problem}

Curran and Locke \cite{CurranLocke} have characterized the $C_4$-face-magic projective labelings on the $4 \times 4$ projective grid graph $\mathcal{P}_{4,4}$.
They show that there are 144 $C_4$-face-magic projective labelings on $\mathcal{P}_{4,4}$ up to symmetries on the projective plane.

\end{document}